\documentclass{article}
\usepackage{latexsym,amsfonts,amsmath,amsthm,amssymb,makeidx}
\usepackage[title]{appendix}
\usepackage{CJK,CJKnumb,CJKulem,times,ifthen,mathrsfs,latexsym,amsfonts, color}
\usepackage{amsmath,amsthm,makeidx,fontenc,amssymb,bm,graphicx,psfrag,listings, curves,extarrows}
\usepackage[driver=pdftex,lmargin=3.2cm,rmargin=3.2cm,heightrounded=true,centering]{geometry}
\usepackage{authblk}
\let\oldbibliography\thebibliography
\renewcommand{\thebibliography}[1]{
\oldbibliography{#1}
\setlength{\itemsep}{0pt}
}

\newcommand{\be}{\begin{equation}}
\newcommand{\ee}{\end{equation}}
\newcommand{\bea}{\begin{eqnarray}}
\newcommand{\eea}{\end{eqnarray}}
\newcommand{\bes}{\begin{eqnarray*}}
\newcommand{\ees}{\end{eqnarray*}}

\newcommand{\nn}{\nonumber}

\def\>{\rangle}
\newcommand{\lb}{\label}

\newcommand{\hb}{\vrule height0.18cm width0.14cm $\,$}

\numberwithin{equation}{section}
\begin{document}
\title{{\bf Morse index and Maslov-type index of the discrete Hamiltonian system }
\footnote{G. Zhu is partially supported by NSF of China
(11871356).}}

\author[1]{Gaosheng Zhu\thanks{gaozsc@163.com}}
\affil[1]{School of Mathematics, Tianjin University, Tianjin 300072, P. R. China}
\date{}
\maketitle




\begin{abstract}
{\it  In this paper, we give the definition of Maslov-type index of the discrete Hamiltonian system, and obtain the relation of Morse index and Maslov-type index of the discrete Hamiltonian system which is a generalization of the case $\omega=1$ in \cite{RoS1}, \cite{RoS2} and \cite{Maz1} to the case $\omega \in {\bf U}$ via direct method which is different from that of \cite{RoS1}, \cite{RoS2} and \cite{Maz1}. Moreover the well-posedness of the splitting numbers $\mathcal {S}_{h,\omega}^{\pm}$ given by (\ref{1.61}) is proven, thus the index iteration theories in \cite{Bot1} and \cite{Lon4} are also valid for the discrete case.}
\end{abstract}



\renewcommand{\theequation}{\thesection.\arabic{equation}}
\renewcommand{\thefigure}{\thesection.\arabic{figure}}

\setcounter{equation}{0}
\section{Introduction and main results}
The main purpose of this paper is to investigate the relation of Morse index and Maslov-type index of the discrete linear Hamiltonian system. The continuous and discrete Hamiltonian system have been studied extensively by many authors, such as the existence and multiplicity of periodic solutions of Hamiltonian system etc. Since what we discussed in this paper is the relation of Morse index and Maslov-type index, we only need to deal with the linear Hamiltonian system case here.

On one hand, for continuous case, we consider the boundary problem of the linear continuous Hamiltonian system
\be
        \dot {z}=JB(t)z(t) \label{1.1}
\ee
with $z(1)=\omega z(0)$, where $z(t)=(x^{T}(t),y^{T}(t))^{T}$ and $x(t), y(t) \in {\bf C}^{m},$
\be
J=\left(
              \begin{array}{cc}
                0 & -I_{m} \\
                I_{m} & 0 \\
              \end{array}
            \right)\lb{1.2}
\ee
with $I_{m}$ is the identity matrix in ${\bf R}^m$ and
\be
B(t)=\left(
              \begin{array}{cc}
                A(t) & C(t) \\
                C^{T}(t) & D(t) \\
              \end{array}
            \right),\lb{1.3}
\ee
here $B(t)^{T}=B(t)$ and $A(t), C(t)$ $ D(t)$ are $m\times m$ matrices in ${\bf R}^{m}$, $t \in {\bf R}/{\bf Z}, \omega \in {\bf U}=\{\omega \in {\bf C}||\omega|=1\}$.\\
The action functional corresponding to (\ref{1.1})
\bea
\tilde{\mathcal{A}}(z)=\frac{1}{2}\int_{0}^{1}[(-J{\dot z}(t),z(t))-(B(t)z(t),z(t))]dt,\lb{1.65}
\eea
on $W^{1,2}({\bf R}/{\bf Z}, {\bf C}^{2m})=\{z|\int _{0}^{1}(|z|^2+|\dot z|^{2})dt< +\infty, z(1)=\omega z(0), \omega \in {\bf U}\}$. Since the functional $\tilde{\mathcal{A}}$ is strongly indefinite, Morse index is infinite. Whereas Morse index of the reduced functional is finite, the relation of Morse index of the reduced functional corresponding to (\ref{1.1}) and Maslov-type index of continuous Hamiltonian system was established in \cite{Abb1}, \cite{CoZ1}, \cite{Dui1}, \cite{Lon1}, \cite{Lon2}, \cite{Lon3} and \cite{Vit1}.

On the other hand, it is well known that the discrete Hamiltonian system has many applications in computation mathematics,  but also it can be applied to abstract mathematics (\cite{RoS1}, \cite{RoS2}, \cite{Maz1} etc). There are some results on the relation of Morse index and Maslov-type index for the discrete Hamiltonian system, which can be employed to prove the Conley conjecture for the torus case such as in \cite{Maz1}. In \cite{RoS1} and \cite{RoS2}, one result similar to that of continuous Hamiltonian system (\ref{1.1}) was obtained for  $\omega=1$ nondegenerate case, which was generalized to $\omega=1$ degenerate case in \cite{Maz1}. This paper extends the relevant results of \cite{RoS1}, \cite{RoS2} and \cite{Maz1} to $\omega \in {\bf U}$ degenerate case via direct method which is different from that of \cite{RoS1}, \cite{RoS2} and \cite{Maz1}.

To be precise, consider the boundary of the discrete linear Hamiltonian system
\be
               \frac{z_{n+1}-z_{n}}{h}=JB_{n}\tilde {z}_{n} \lb{1.66}
               \ee
with $\tilde z_{n+N}=\omega \tilde z_{n}$, where $0 < h=\frac{1}{N}\in {\bf R}$,  $z_{n}=(x^{T}_{n},y^{T}_{n})^{T}$, $\tilde {z}_{n}=(x^{T}_{n+1},y^{T}_{n})^{T}$ with $x_{n}, y_{n} \in {\bf C}^{m}$, $N \in {\bf N}$, and
\be
B_{n}=\left(
              \begin{array}{cc}
                A_{n} & C_{n} \\
                C^{T}_{n} & D_{n} \\
              \end{array}
                        \right)\lb{1.67}
\ee
satisfying
\bea
B_{n}=B_{n+N}, ~~~B^{T}_{n}=B_{n},\lb{1.291}
\eea
where $A_{n}, C_{n}$ and $D_{n}$ are $m\times m$ real matrices.

Next we introduce two spaces
\bea
{W}=\{z=({z}_{1},{z}_{2},\dots,{z}_{N})|{ z}_{n}=\omega{ z}_{n+N}, {z}_{n}=(x_{n}^{T},y_{n}^{T})^{T}, x_{n},y_{n} \in {\bf C}^{m}, \forall n \in {\bf N}, \omega\in {\bf U}\}\lb{1.292}
\eea
and
\bea
\tilde{W}=\{\tilde {z}=({\tilde z}_{1},{\tilde z}_{2},\dots,{\tilde z}_{N})|{\tilde z}_{n}=\omega{\tilde z}_{n+N}, {\tilde z}_{n}=(x_{n+1}^{T},y_{n}^{T})^{T}, x_{n},y_{n} \in {\bf C}^{m}, \forall n \in {\bf N}, \omega\in {\bf U}\},\lb{1.46}
\eea
It is obvious that there is a linear isomorphism between $W$ and $\tilde W$, hence we can introduce the following action functional
\bea
\tilde{\mathcal{A}}_{h,\omega}(\tilde {z})=\frac{1}{2}\Sigma_{n=1}^{N}({-J\frac{{z}_{n+1}-z_{n}}{h},{\tilde z}_{n}})-(B_{n}\tilde{z}_{n}, {\tilde z}_{n}),\lb{1.294}
\eea
where $B_{n}$ is given by (\ref{1.67}) and satisfies the condition (\ref{1.291}).

Correspondingly, the solution of the boundary problem (\ref{1.66}) is one to one correspondence to the the critical point of $\tilde{\mathcal{A}}_{h,\omega}$.\\
The Hessian of $\tilde{\mathcal{A}}_{h,\omega}$ at $\tilde {z}=0$ is
\bea
\mathcal{A}_{h,\omega}(\tilde {z},\tilde {w})=\Sigma_{n=1}^{N}({-J\frac{{z}_{n+1}-z_{n}}{h},\tilde{w}_{n}})-(B_{n}\tilde{z}_{n} \tilde{w}_{n}) \lb{1.47}
\eea
on $\tilde{W}$. Denote by $m_{h,\omega}^{-}, m_{h,\omega}^{0}, m_{h,\omega}^{+}$ the negative, null, positive Morse indices of $\tilde{\mathcal{A}}_{h,\omega}$, that is, the dimension of negative, zero, positive invariant subspace of ${\mathcal{A}}_{h,\omega}$ respectively. Moreover, we can associate index functions $i_{\omega}(\gamma_{d,h,N}),\nu_{\omega}(\gamma_{d,h,N})$ to the discrete Hamiltonian system (\ref{1.66}), note that as $N \rightarrow \infty$, there hold $m_{h,\omega}^{-}, m_{h,\omega}^{+} \rightarrow \infty$, but $Sign \mathcal{A}_{h,\omega}$ is finite. To be precise, we obtain the following\\

{\it{\bf Theorem 1.1.} For any $\omega \in {\bf U}$, there hold
\begin{eqnarray}
m_{h,\omega}^{-}&=&mN+i_{\omega}(\gamma_{d,h,N}),\nn\\
m_{h,\omega}^{0}&=&\nu_{\omega}(\gamma_{d,h,N}),\nn\\
m_{h,\omega}^{+}&=& mN-i_{\omega}(\gamma_{d,h,N})-\nu_{\omega}(\gamma_{d,h,N}).\lb{1.64}
\end{eqnarray}
where $h=\frac{1}{N}, N\in {\bf N}$ is sufficiently large.}

As consequences of Theorem 1.1, we have two corollaries.\\

{\it{\bf Corollary 1.2.} (cf. Theorem 4.1 of \cite{RoS1} or Theorem 6.4 of \cite{RoS2} for the case $\omega=1$ ) As situations above,\\
i). if $\nu_{\omega}(\gamma_{d,h,N})=0$, we have
\bea
Sign \mathcal{A}_{h,\omega}=2i_{\omega}(\gamma_{d,h,N}).\lb{1.68}
\eea
ii). if $B_{n}=B(\frac{n}{N})$ and $\nu_{\omega}(\gamma)=0$, we have
\bea
\nu_{\omega}(\gamma_{d,h,N})&=&0, \nn\\
i_{\omega}(\gamma_{d,h,N})&=&i_{\omega}(\gamma), \nn\\
Sign \mathcal{A}_{h,\omega}&=&2i_{\omega}(\gamma).\lb{1.293}
\eea
where $B(\cdot)$ and $\gamma$ are the coefficient matrix and fundamental solution of (\ref{1.1}) respectively, $h=\frac{1}{N}, N\in {\bf N}$ is sufficiently large, $Sign \mathcal{A}_{h,\omega}=m_{h,\omega}^{-}-m_{h,\omega}^{+}$.}\\

The other corollary is concerned with the splitting numbers $\mathcal{S}_{h,\omega}^{\pm}$ and $S_{M}^{\pm}(\omega)$, which are given by (\ref{1.61}) and (\ref{1.297}) respectively. To be Precise, we obtain\\

{\it{\bf Corollary 1.3.} Under the same conditions as Corollary 1.2, then\\
i). $\mathcal{S}_{h,\omega}^{\pm}=S_{\gamma_{c,h,N}(1)}^{\pm}(\omega)$.\\
ii). if $B_{n}=B(\frac{n}{N})$ and $\nu_{\omega}(\gamma)=0$, we have
$\mathcal{S}_{h,\omega}^{\pm}=S_{\gamma(1)}^{\pm}(\omega)$.\\
Where $B(\cdot)$ and $\gamma$ are the coefficient matrix and fundamental solution of (\ref{1.1}) respectively, $h=\frac{1}{N}, N\in {\bf N}$ is sufficiently large.}\\

{\it {\bf Remark 1.4.} i). The second item of Corollary 1.2 is concerned with the relation of Maslov-type index the continuous Hamiltonian system and the corresponding discrete Hamiltonian system, especially shows that $\nu_{\omega}(\gamma_{d,h,N}),i_{\omega}(\gamma_{d,h,N})$ and $Sign \mathcal{A}_{h,\omega}$ possess some invariant property independent of the choice of sufficiently small $h$.\\
ii). Based on Theorem 1.1 and Lemma 4.6, the splitting numbers $\mathcal{S}_{h,\omega}^{\pm}$ are only dependent of $\gamma_{c,h,N}(1)$, but independent of the path $\gamma_{c,h,N}$,  then they possess same properties as those for the continuous case, hence are well-defined. Moreover Corollary 1.3 implies that the index iteration theories in \cite{Bot1} and \cite{Lon4} are also valid for the discrete Hamiltonian system case. Especially, in the case ii) of Corollary 1.3, $\mathcal{S}_{h,\omega}^{\pm}$ are independent of $h$.}\\

This paper is motivated by \cite{Lon1}, \cite{Lon2}, \cite{Lon3}, \cite{LoZ1}, \cite{Lon4}, \cite{Maz1}, \cite{RoS1} and \cite{RoS2}. The paper is organized as follows. We provide a detailed exposition of the equivalence between the symplecticity of $S_{n}$ given by (\ref{1.22}) and the symmetry of $B_{n}$ given by (\ref{1.5}), then introduce the definition of Maslov-type index of the discrete Hamiltonian system and prove that it is well defined in the second section. Section 3 gives the presentation of the eigenvalues and eigenvectors of discrete Hamiltonian system. We reviews without proofs some relevant materials on Maslov index and verify some properties of indices for discrete Hamiltonian system in the fourth section. Section 5 is devoted to the construction of the homotopy of discrete symplectic paths which is used to reduce the proof of Theorem 1.1 for general discrete Hamiltonian system to the proof for the discrete Hamiltonian with coefficient matrix ${\hat B}_{j,n}$ which corresponds to ${\hat B}_{j}(t)$ associated to the standard path ${\hat \beta}_{j}(t)$. Section 6 is intended to prove Theorem 1.1 and Corollary 1.3. The estimates of the errors between the continuous Hamiltonian system and the corresponding discrete Hamiltonian system are given in Section 7.

In this paper, let ${\bf N}, {\bf Z}$, {\bf R} and ${\bf C}$ denote the sets of natural integers, integers, real and complex numbers respectively.

\setcounter{equation}{0}
\section{Definition of Maslov-type index of discrete Hamiltonian system}
In this section we consider the following discrete system
\be
\frac{z_{n+1}-z_{n}}{h}=JB_{n}\tilde {z}_{n},\lb{1.4}
\ee
where $0 < h\in {\bf R}$,  $z_{n}=(x^{T}_{n},y^{T}_{n})^{T}$, $\tilde {z}_{n}=(x^{T}_{n+1},y^{T}_{n})^{T}$ with $x_{n}, y_{n} \in {\bf R}^{m}$, and
\be
B_{n}=\left(
              \begin{array}{cc}
                A_{n} & C_{n} \\
                F_{n} & D_{n} \\
              \end{array}
            \right)\lb{1.5}
\ee
here $A_{n}, C_{n}, F_{n}$ and $D_{n}$ are $m\times m$ real matrices, $n \in {\bf N}$. Since $B_{n}$ are real matrices, the Maslov-type index is independent on the fact that the domain of entries of symplectic matrices is either {\bf R} or ${\bf C}$, thus we can assume that $x_{n}, y_{n} \in {\bf R}^{m}$.

\subsection {The relation between the symplecticity of $S_{n}$ and the symmetry of $B_{n}$}
Since the following computations of this section for the general $h$ are similar to those for $h=1$, for simplicity of notation, we take $h=1$.  We can write (\ref {1.4}) in the form
\begin{eqnarray}
\left(
              \begin{array}{cc}
                x_{n+1}-x_{n}  \\
                y_{n+1}-y_{n} \\
              \end{array}
            \right)
           &=&\left(
              \begin{array}{cc}
                0 & -I \\
                I & 0 \\
              \end{array}
            \right)
            \left(
              \begin{array}{cc}
                A_{n} & C_{n} \\
                F_{n} & D_{n} \\
              \end{array}
            \right)
            \left(
              \begin{array}{cc}
                x_{n+1}  \\
                y_{n} \\
              \end{array}
            \right)\nn\\
            &=&\left(
            \begin{array}{cc}
            -F_{n}x_{n+1}-D_{n}y_{n}\\
            A_{n}x_{n+1}+C_{n}y_{n}\\
            \end{array}
            \right)\lb{1.6},
\end{eqnarray}
That is
\be
\left\{\begin{array}{ll}
               x_{n+1}-x_{n}=-F_{n} x_{n+1}-D_{n}y_{n},~~~~~~~~~~~~~~~~~~~~~~~~~~~~~~~~~~~~ \\
               y_{n+1}-y_{n}=A_{n} x_{n+1}+C_{n}y_{n}. ~~~~~~~~~~~~~~~~~~~~~~~~~~~~~~~~~~~~~~
               \end{array}
\right. \lb{1.7}\ee
Then we have
\be
\left\{\begin{array}{ll}
               y_{n+1}=A_{n}x_{n+1}+(I+C_{n})y_{n},~~~~~~~~~~~~~~~~~~~~~~~~~~~~~~~~~~~~ \\
               x_{n+1}=(I+F_{n})^{-1}(x_{n}-D_{n}y_{n}).~~~~~~~~~~~~~~~~~~~~~~~~~~~~~~~ ~~~
               \end{array}
\right.\lb{1.8}
\ee
Substituting the second line into the first line yields that
\begin{eqnarray}
y_{n+1}&=&A_{n}[(I+F_{n})^{-1}(x_{n}-D_{n}y_{n})]+(I+C_{n})y_{n}\nn\\
       &=&A_{n}(I+F_{n})^{-1}x_{n}+[-A_{n}(I+F_{n})^{-1}D_{n}+(I+C_{n})]y_{n}, \lb{1.9}
\end{eqnarray}
thus there holds
\begin{eqnarray}
\left(
              \begin{array}{cc}
                x_{n+1}  \\
                y_{n+1} \\
              \end{array}
            \right)
            =\left(
              \begin{array}{cc}
                (I+F_{n})^{-1} & -(I+F_{n})^{-1}D_{n} \\
                A_{n}(I+F_{n})^{-1} & -A_{n}(I+F_{n})^{-1}D_{n}+I+C_{n} \\
              \end{array}
            \right)
            \left(
              \begin{array}{cc}
                x_{n} \\
                y_{n} \\
              \end{array}
            \right).\lb{1.10}
\end{eqnarray}
Write (\ref{1.10}) as
\bea
z_{n+1}=S_{n}z_{n},\lb{1.131}
\eea
where
\begin{eqnarray}
S_{n}:=\left(
               \begin{array}{cc}
               S_{n,11} & S_{n,12}\\
               S_{n,21} & S_{n,22}\\
               \end{array}
             \right)
       =\left(
              \begin{array}{cc}
                (I+F_{n})^{-1} & -(I+F_{n})^{-1}D_{n} \\
                A_{n}(I+F_{n})^{-1} & -A_{n}(I+F_{n})^{-1}D_{n}+I+C_{n} \\
              \end{array}
            \right),\lb{1.22}
\end{eqnarray}
and $S_{n,11}, S_{n,12}, S_{n,21}, S_{n,22}$ are $m\times m$ matrices.

Moreover, denote
\bea
z_{n+1}^{l}=S_{n}z_{n}^{l}, \lb{1.132}
\eea
by
\be
\gamma_{n+1}=S_{n}\gamma_{n},\lb{1.40}
\ee
where $\gamma_{n}=(z_{n}^{1}, z_{n}^{2}, \cdots z_{n}^{2m}), \gamma_{n+1}=(z_{n+1}^{1}, z_{n+1}^{2}, \cdots z_{n+1}^{2m})$ are $2m\times 2m$ matrices, $z_{n}^{l}, z_{n+1}^{l} \in {\bf R}^{2m}$, $l=1,2,\cdots,2m$. Then we have\\

{\it {\bf Lemma 2.1.} If $S_{n}$ is defined by (\ref{1.22}), $S_{n,11}$ is invertible, then the matrix $S_{n}$ is symplectic if and only if $A_{n}^{T}=A_{n}, F_{n}=C_{n}^{T}, D_{n}^{T}=D_{n}$. }

{\bf Proof.} From the definition of $S_{n}$, we have $(I+F_{n})^{-1}=S_{n,11}$, it follows that
\be
F_{n}=S_{n,11}^{-1}-I.\lb{1.14}
\ee
Also, since $-(I+F_{n})^{-1}D_{n}=S_{n,12}$, there holds
\be
D_{n}=-S_{n,11}^{-1}S_{n,12}.\lb{1.15}
\ee
Meanwhile, (\ref{1.22}) yields $A_{n}(I+F_{n})^{-1}=S_{n,21}$, it gives
\be
A_{n}=S_{n,21}S_{n,11}^{-1}.\lb{1.16}
\ee
Finally, according to $-A_{n}(I+F_{n})^{-1}D_{n}+I+C_{n}=S_{n,22}$, we see that
\be
C_{n}=S_{n,22}-I+A_{n}S_{n,11}D_{n}.\lb{1.17}
\ee

By Lemma 1.1.2 of \cite{Lon4}, the fact that $S_{n}$ is symplectic is equivalent to the following three equations hold
\begin{eqnarray}
& &S_{n,11}^{T}S_{n,22}-S_{21}^{T}S_{n,12}=I, \lb{1.11}\\
& &S_{n,11}^{T}S_{n,21}=S_{n,21}^{T}S_{n,11}, \lb{1.12}\\
& &S_{n,12}^{T}S_{n,22}=S_{n,22}^{T}S_{n,12}, \lb{1.13}
\end{eqnarray}
therefore, we only need to show that (\ref{1.11}), (\ref{1.12}),(\ref{1.13}) are true if and only if $A_{n}^{T}=A_{n}, F_{n}=C_{n}^{T}, D_{n}^{T}=D_{n}$, The proof is divided into three steps.

{\bf Step 1.} The equivalence of (\ref{1.12}) and $A_{n}^{T}=A_{n}$.

From (\ref{1.12}), we have $(S_{n,11}^{T})^{-1}S_{n,21}^{T}=S_{n,21}S_{n,11}^{-1}$, that is $(S_{n,11}^{-1})^{T}S_{n,21}^{T}=S_{n,21}S_{n,11}^{-1}$, moreover there holds $(S_{n,21}S_{n,11}^{-1})^{T}=S_{n,21}S_{n,11}^{-1}$, together (\ref{1.16}), thus $A_{n}^{T}=A_{n}$. Since the process above can be reversible, (\ref{1.12}) is equivalent to $A_{n}^{T}=A_{n}$.

{\bf Step 2.} (\ref{1.11}), (\ref{1.12}) hold if and only if $F_{n}^{T}=C_{n}$.

By (\ref{1.11}), there holds
\begin{eqnarray}
I&=&S_{n,11}^{T}S_{n,22}-S_{n,21}^{T}S_{n,12}\nn\\
 &=&S_{n,11}^{T}S_{n,22}-S_{n,21}^{T}S_{n,11}S_{n,11}^{-1}S_{n,12}\nn\\
 &=&S_{n,11}^{T}S_{n,22}-S_{n,11}^{T}S_{n,21}S_{n,11}^{-1}S_{n,12},\lb{1.18}
\end{eqnarray}
where we use (\ref{1.12}) in the third equality. Combining (\ref{1.15}) with  (\ref{1.16}),  from (\ref{1.18}), thus we get that
\begin{eqnarray}
(S_{n,11}^{T})^{-1}&=&S_{n,22}-S_{n,21}S_{11}^{-1}S_{n,12}\nn\\
                 &=&S_{n,22}+S_{n,21}S_{n,11}^{-1}S_{n,11}(-S_{n,11}^{-1})S_{n,12}\nn\\
                 &=&S_{n,22}+A_{n}S_{n,11}D_{n}.\lb{1.19}
\end{eqnarray}
Meanwhile, it follows that
\begin{eqnarray}
F_{n}^{T}&=&(S_{n,11}^{-1}-I)^{T}\nn\\
         &=&(S_{n,11}^{T})^{-1}-I.\lb{1.20}
\end{eqnarray}
from (\ref{1.14}).

Finally, by  (\ref{1.17}), (\ref{1.19}) and (\ref{1.20}), we have $F_{n}^{T}=C_{n}$. Since the process above can be reversible, (\ref{1.11}) is equivalent to $F_{n}^{T}=C_{n}$.

{\bf Step 3.} It remains to prove that (\ref{1.13}) is equivalent to $D_{n}^{T}=D_{n}$.

From (\ref{1.22}) and (\ref{1.13}), it leads to
\begin{eqnarray}
& &[-(I+F_{n})^{-1}D_{n}]^{T}[-A_{n}(I+F_{n})^{-1}D_{n}+(I+C_{n})]\nn\\
                          \hspace*{0.2cm}&&=[-A_{n}(I+F_{n})^{-1}D_{n}+(I+C_{n})]^{T}[-(I+F_{n})^{-1}D_{n}].\lb{1.21}
\end{eqnarray}
By simple computations, we see that
\begin{eqnarray}
& &D_{n}^{T}[(I+F_{n})^{-1}]^{T}[A_{n}(I+F_{n})^{-1}D_{n}]-D_{n}^{T}[(I+F_{n})^{-1}]^{T}(I+C_{n})\nn\\
&&=D_{n}^{T}[(I+F_{n})^{-1}]^{T}[A_{n}^{T}(I+F_{n})^{-1}D_{n}]-(I+C_{n})^{T}(I+F_{n})^{-1}D_{n}\nn\\
&&=D_{n}^{T}[(I+F_{n})^{-1}]^{T}[A_{n}(I+F_{n})^{-1}D_{n}]-(I+C_{n})^{T}(I+F_{n})^{-1}D_{n},\nn
\end{eqnarray}
where $A_{n}^{T}=A_{n}$ is used in the second equality. Therefore, there holds
\begin{eqnarray}
D_{n}^{T}[(I+F_{n})^{-1}]^{T}(I+C_{n})=(I+C_{n})^{T}(I+F_{n})^{-1}D_{n},\nn
\end{eqnarray}
together with $F_{n}^{T}=C_{n}$ in step 2, we have $D_{n}^{T}=D_{n}$.\hfill\hb

Set
\begin{eqnarray}
X_{n}=\left(
               \begin{array}{cc}
               X_{n,11} & X_{n,12}\\
               X_{n,21} & X_{n,22}\\
               \end{array}
             \right),\lb{1.50}
\end{eqnarray}
where $X_{n,11}, X_{n,12}, X_{n,21}, X_{n,22}$ are $m\times m$ matrices. The following statement holds\\

{\it {\bf Lemma 2.2.} Suppose that $X_{n}$ is close enough to $I_{2m}$, then\\
{\bf i).} The fact that $X_{n}$ is symplectic if and only if there exist sufficiently close to zero $m\times m$ matrices ${\tilde A }_{n}, {\tilde C }_{n}, {\tilde D }_{n}$ with ${\tilde A }_{n}^{T}={\tilde A }_{n}, {\tilde D }_{n}^{T}={\tilde D}_{n}$ such that $X_{n}$ has the form
\begin{eqnarray}
X_{n}
            =\left(
              \begin{array}{cc}
                (I+{\tilde C }_{n}^{T})^{-1} & -(I+{\tilde C }_{n}^{T})^{-1}{\tilde D }_{n} \\
                {\tilde A }_{n}(I+{\tilde C }_{n}^{T})^{-1} & -{\tilde A }_{n}(I+{\tilde C }_{n}^{T})^{-1}{\tilde D }_{n}+I+{\tilde C }_{n} \\
              \end{array}
            \right).\lb{1.51}
\end{eqnarray}
{\bf ii)}. There is one to one correspondence between the symplectic matrices $X_{n}$ and the symmetric matrices
\be
{\tilde B }_{n}=\left(
              \begin{array}{cc}
                {\tilde A }_{n} & {\tilde C }_{n} \\
                {\tilde C }_{n}^{T} & {\tilde D }_{n} \\
              \end{array}
            \right)\nn
\ee
as $B_{n}$ in (\ref{1.4}). }

{\bf Proof. i).} Sufficiency is followed from Lemma 2.1, we only need to prove the necessity.\\
Since $X_{n}$ is sufficiently close to $I_{2m}$, firstly we can define $C_{n}$ sufficiently small as
\bea
{\tilde C }_{n}=(X_{n,11}^{-1}-I)^{T},\lb{1.194}
\eea
then define ${\tilde A }_{n}$ small enough satisfying
\bea
{\tilde A }_{n}=X_{n,21}(I+{\tilde C }_{n}^{T}),\lb{1.195}
\eea
finally define ${\tilde D }_{n}$ sufficiently small to be
\bea
{\tilde D }_{n}=-(I+{\tilde C }_{n}^{T})X_{n,12},\lb{1.196}
\eea
thus the proof is reduced to verifying that $X_{n,22}=-{\tilde A }_{n}(I+{\tilde C }_{n}^{T})^{-1}{\tilde D }_{n}+I+{\tilde C }_{n},  {\tilde A }_{n}^{T}={\tilde A }_{n}$ and ${\tilde D }_{n}^{T}={\tilde D }_{n}$.

By Lemma 1.1.2 of \cite{Lon4}, the fact that $X_{n}$ is symplectic is equivalent to the following three equations hold
\begin{eqnarray}
& &X_{n,11}^{T}X_{n,22}-X_{n,21}^{T}X_{n,12}=I, \lb{1.52}\\
& &X_{n,11}^{T}X_{n,21}=X_{n,21}^{T}X_{n,11}, \lb{1.53}\\
& &X_{n,12}^{T}X_{n,22}=X_{n,22}^{T}X_{n,12}, \lb{1.54}
\end{eqnarray}
The proof is divided into three steps.

{\bf Step 1.} ${\tilde A }_{n}^{T}={\tilde A }_{n}$.\\
Substituting $X_{n,11}, X_{n,21}$ into (\ref{1.53}), we obtain
\bea
[(I+{\tilde C }_{n}^{T})^{-1}]^{T}{\tilde A }_{n}(I+{\tilde C }_{n}^{T})^{-1}=[{\tilde A }_{n}(I+{\tilde C }_{n}^{T})^{-1}]^{T}(I+{\tilde C }_{n}^{T})^{-1},\nn
\eea
that is
\bea
(I+{\tilde C }_{n})^{-1}{\tilde A }_{n}(I+{\tilde C }_{n}^{T})^{-1}=(I+{\tilde C }_{n})^{-1}{{\tilde A }_{n}}^{T}(I+{\tilde C }_{n}^{T})^{-1},\nn
\eea
hence ${\tilde A }_{n}^{T}={\tilde A }_{n}$ holds.

{\bf Step 2.} $X_{n,22}=-{\tilde A }_{n}(I+{\tilde C }_{n}^{T})^{-1}{\tilde D }_{n}+I+{\tilde C }_{n}$.\\
By (\ref{1.52}), we see that
\bea
[(I+{\tilde C }_{n}^{T})^{-1}]^{T}X_{22}-[{\tilde A }_{n}(I+{\tilde C }_{n}^{T})^{-1}]^{T}[-(I+{\tilde C }_{n}^{T})^{-1}{\tilde D }_{n}]=I,\nn
\eea
This gives
\bea
X_{n,22}&=&(I+{\tilde C }_{n})\{I-[{\tilde A }_{n}(I+{\tilde C }_{n}^{T})^{-1}]^{T}(I+{\tilde C }_{n}^{T})^{-1}{\tilde D }_{n}\}\nn\\
      &=& (I+{\tilde C }_{n})\{I-(I+{\tilde C }_{n})^{-1}{\tilde A }_{n}^{T}(I+{\tilde C }_{n}^{T})^{-1}{\tilde D }_{n}\}\nn\\
      &=& (I+{\tilde C }_{n})-{\tilde A }_{n}^{T}(I+{\tilde C }_{n}^{T})^{-1}{\tilde D }_{n}\nn\\
      &=& (I+{\tilde C }_{n})-{\tilde A }_{n}(I+{\tilde C }_{n}^{T})^{-1}{\tilde D }_{n},\nn
\eea
we use ${\tilde A }_{n}^{T}={\tilde A }_{n}$ in the last equality.

{\bf Step 3.} ${\tilde D }_{n}^{T}={\tilde D }_{n}$.\\
Substituting $X_{n,12}, X_{n,22}$ into (\ref{1.54}) yields that
\bea
[(I+{\tilde C }_{n}^{T})^{-1}{\tilde D }_{n}]^{T}[-{\tilde A }_{n}(I+{\tilde C }_{n}^{T})^{-1}{\tilde D }_{n}+I+{\tilde C }_{n}]=[-{\tilde A }_{n}(I+{\tilde C }_{n}^{T})^{-1}{\tilde D}_{n}+I+{\tilde C }_{n}]^{T}[(I+{\tilde C }_{n}^{T})^{-1}{\tilde D }_{n}].\nn
\eea
By simple computation, we have
\bea
-{\tilde D }_{n}^{T}(I+{\tilde C }_{n})^{-1}{\tilde A }_{n}(I+{\tilde C }_{n}^{T})^{-1}{\tilde D }_{n} + {\tilde D }_{n}^{T}=-{\tilde D }_{n}^{T}(I+{\tilde C }_{n})^{-1}{\tilde A }_{n}^{T}(I+{\tilde C }_{n}^{T})^{-1}{\tilde D }_{n} + {\tilde D }_{n},\nn
\eea
together with ${\tilde A }_{n}^{T}={\tilde A }_{n}$, there holds ${\tilde D }_{n}^{T}={\tilde D }_{n}$.

{\bf ii).} It is followed from the definition of $B_{n}$ and Lemma 2.1.\hfill\hb\\

{\it {\bf Remark 2.3. i)}. If ${\tilde {z}}_{n}$ of the right hand side of (\ref{1.4}) is replaced by $z_{n}$, the corresponding $S_{n}$ is not necessarily a symplectic matrix.\\
{\bf ii).} For the general $h$, the correspondence with (\ref{1.22}) is
\begin{eqnarray}
S_{n,h}:&=&\left(
               \begin{array}{cc}
               S_{n,11,h} & S_{n,12,h}\\
               S_{n,21,h} & S_{n,22,h}\\
               \end{array}
             \right)\nn\\
       &=&
\left(
              \begin{array}{cc}
                (I+hF_{n})^{-1} & -h(I+hF_{n})^{-1}D_{n} \\
                hA_{n}(I+hF_{n})^{-1} & -h^{2}A_{n}(I+hF_{n})^{-1}D_{n}+I+hC_{n} \\
              \end{array}
            \right),\lb{1.37}
            \end{eqnarray}
and
\be
\gamma_{n+1}=S_{n,h}\gamma_{n}.\lb{1.42}
\ee
{\bf iii).} It follows that  $S_{n,h}$ are symplectic matrices if and if only $A_{n}^{T}=A_{n}, F_{n}=C_{n}^{T}, D_{n}^{T}=D_{n}$ from Lemma 2.1.\\
{\bf iv)}. For the general $h$, the correspondence with (\ref{1.194}), (\ref{1.195}), (\ref{1.196}) is respectively
\bea
{h\tilde C }_{n}=(X_{n,11}^{-1}-I)^{T},\lb{1.204}
\eea
\bea
{h\tilde A }_{n}=X_{n,21}(I+h{\tilde C }_{n}^{T}),\lb{1.205}
\eea
\bea
{h\tilde D }_{n}=-(I+h{\tilde C }_{n}^{T})X_{n,12},\lb{1.206}
\eea
{\bf v).} If the discrete Hamiltonian system is
\be
\frac{z_{n+1}-z_{n}}{h}=JB_{r,n}\tilde {z}_{n}\lb{1.211}
\ee
with
\be
B_{r,n}=\left(
              \begin{array}{cc}
                A_{r,n} & C_{r,n} \\
                F_{r,n} & D_{r,n} \\
              \end{array}
            \right),\lb{1.212}
\ee
the counterpart of (\ref{1.37}), (\ref{1.42}) is respectively
\begin{eqnarray}
S_{r,n,h}=
\left(
              \begin{array}{cc}
                (I+hF_{r,n})^{-1} & -h(I+hF_{r,n})^{-1}D_{r,n} \\
                hA_{r,n}(I+hF_{r,n})^{-1} & -h^{2}A_{r,n}(I+hF_{r,n})^{-1}D_{r,n}+I+hC_{r,n} \\
              \end{array}
            \right),\lb{1.213}
            \end{eqnarray}
\be
\gamma_{r,n+1}=S_{r,n,h}\gamma_{r,n},\lb{1.214}
\ee
Meanwhile, if $X_{n}$ in (\ref{1.50}) is replaced by
\begin{eqnarray}
X_{r,n}=\left(
               \begin{array}{cc}
               X_{r,n,11} & X_{r,n,12}\\
               X_{r,n,21} & X_{r,n,22}\\
               \end{array}
             \right),\lb{1.221}
\end{eqnarray}
then (\ref{1.204}), (\ref{1.205}), (\ref{1.206}) is respectively corresponded with
\bea
{h\tilde C }_{r,n}=(X_{r,n,11}^{-1}-I)^{T},\lb{1.208}
\eea
\bea
{h\tilde A }_{r,n}=X_{r,n,21}(I+h{\tilde C }_{r,n}^{T}),\lb{1.219}
\eea
\bea
{h\tilde D }_{r,n}=-(I+h{\tilde C }_{r,n}^{T})X_{r,n,12},\lb{1.220}
\eea}

\subsection{The definition of Maslov-type index of discrete Hamiltonian system}
Based on the discussion of Subsection 2.1, here and in the sequel, we always suppose that $B_{n}=B_{n}^T$. To find a family of continuous symplectic matrices which connect $\gamma_{n}$ to $\gamma_{n+1}$, the basic idea is as follows.

Suppose
\begin{eqnarray}
\left(
              \begin{array}{cc}
                x_{n,s}-x_{n}  \\
                y_{n,s}-y_{n} \\
              \end{array}
            \right)
           &=&s
           \left(
              \begin{array}{cc}
                0 & -I \\
                I & 0 \\
              \end{array}
            \right)
            \left(
              \begin{array}{cc}
                A_{n} & C_{n} \\
                C_{n}^{T} & D_{n} \\
              \end{array}
            \right)
            \left(
              \begin{array}{cc}
                x_{n,s}  \\
                y_{n} \\
              \end{array}
            \right)\lb{1.33},
\end{eqnarray}
where
\bea
x_{n,s}=\left\{\begin{array}{ll}
x_{n}, ~~~s=0. \\
x_{n+1}, ~s=1.
\end{array}
\right.\nn\eea
and
\bea
y_{n,s}=\left\{\begin{array}{ll}
 y_{n}, ~~~s=0. \\
y_{n+1},~ s=1.\end{array}
\right.\nn\eea
 with $0\le s\le 1.$

Also we can obtain an analogue to (\ref{1.10}) from (\ref{1.33}),
\begin{eqnarray}
\left(
              \begin{array}{cc}
                x_{n,s}  \\
                y_{n,s} \\
              \end{array}
            \right)
            =S_{n,s}
            \left(
              \begin{array}{cc}
                x_{n} \\
                y_{n} \\
              \end{array}
            \right),\lb{1.34}
\end{eqnarray}
where
\begin{eqnarray}
S_{n,s}=
\left(
\begin{array}{cc}
(I+sC_{n}^{T})^{-1} & -s(I+sC_{n}^{T})^{-1}D_{n}\\
sA_{n}(I+sC_{n}^{T})^{-1} & -s^{2}A_{n}(I+sC_{n}^{T})^{-1}D_{n}+(I+sC_{n})\\
\end{array}
\right)\lb{1.35}
\end{eqnarray}
with $0\le s\le 1$.

Moreover note that $S_{n,s}$ is symplectic and
\bea
S_{n,s}=\left\{\begin{array}{ll}
I_{2m}, \hspace*{0.2cm}s=0. \\
S_{n},\hspace*{0.3cm} s=1.
\end{array}
\right.\nn\eea
where $S_{n}$ is given by (\ref{1.22}) with $F_{n}=C_{n}^{T}$. Thus we can introduce the following\\

{\it {\bf Lemma 2.4.} Suppose that $A_{n}^{T}=A_{n}, F_{n}=C_{n}^{T}, D_{n}^{T}=D_{n}$ and $\gamma_{n+1}=S_{n}\gamma_{n}$, where $\gamma_{n}$ is a symplectic matrix, $S_{n}$ is given by (\ref{1.22}), then a family of continuous symplectic matrices path $\gamma_{n,s}$ with two end points $\gamma_{n}$ and $\gamma_{n+1}$ is defined by
\begin{eqnarray}
\gamma_{n,s}=
S_{n,s}\gamma_{n},\lb{1.38}
\end{eqnarray}
where $S_{n,s}$ is given by (\ref{1.35}), $0\le s\le 1$.}\\

{\it {\bf Remark 2.5.} {\bf i).} For the general $h$, the correspondence with (\ref{1.35}) is
\begin{eqnarray}
S_{n,s,h}=
\left(
\begin{array}{cc}
(I+shC_{n}^{T})^{-1} & -sh(I+shC_{n}^{T})^{-1}D_{n}\\
shA_{n}(I+shC_{n}^{T})^{-1} & -s^{2}h^{2}A_{n}(I+shC_{n}^{T})^{-1}D_{n}+(I+shC_{n})\\
\end{array}
\right).\lb{1.41}
\end{eqnarray}
Note that when $h$ is small enough, $(I+shC_{n}^{T})^{-1}$  make sense  and $S_{n,s,h}$ are symplectic matrices for $0\le s\le 1$. At the same time, we also have
\bea
S_{n,s,h}=\left\{\begin{array}{ll}
I_{2m}, \hspace*{0.2cm}s=0. \\
S_{n,h},\hspace*{0.3cm} s=1.
\end{array}\right.\lb{1.43}
\eea
where $S_{n,h}$ is given by (\ref{1.37}) with $F_{n}=C_{n}^{T}$.\\
Moreover, (\ref{1.38}) in Lemma 2.4 should be replaced by
\begin{eqnarray}
\gamma_{n,s,h}=
S_{n,s,h}\gamma_{n},\lb{1.39}
\end{eqnarray}
where $0\le s\le 1$. By (\ref{1.43}) and (\ref{1.39}), there holds $\gamma_{n,0,h}=\gamma_{n}, \gamma_{n,1,h}=\gamma_{n+1}=\gamma_{n+1,0,h}$ and $ \gamma_{n+1,1,h}=\gamma_{n+2}.$\\
{\bf ii).} It follows that  $S_{n,s,h}$ are symplectic matrices if and if only $A_{n}^{T}=A_{n}, D_{n}^{T}=D_{n}$ from Lemma 2.1.\\
{\bf iii).} If the discrete Hamiltonian system is given by (\ref{1.211}) and (\ref{1.212}) with $F_{r,n}=C_{r,n}^{T}$, then the counterpart of (\ref{1.41}), (\ref{1.43}), (\ref{1.39}) is respectively
\begin{eqnarray}
S_{r,n,s,h}=
\left(
\begin{array}{cc}
(I+shC_{r,n}^{T})^{-1} & -sh(I+shC_{r,n}^{T})^{-1}D_{n}\\
shA_{r,n}(I+shC_{r,n}^{T})^{-1} & -s^{2}h^{2}A_{r,n}(I+shC_{r,n}^{T})^{-1}D_{r,n}+(I+shC_{r,n})\\
\end{array}
\right),\lb{1.215}
\end{eqnarray}
\bea
S_{r,n,s,h}=\left\{\begin{array}{ll}
I_{2m}, \hspace*{0.2cm}s=0. \\
S_{r,n,h},\hspace*{0.3cm} s=1.
\end{array}
\right.\lb{1.216}
\eea
\begin{eqnarray}
\gamma_{r,n,s,h}=
S_{r,n,s,h}\gamma_{r,n},\lb{1.217}
\end{eqnarray}
with $\gamma_{r,n,0,h}=\gamma_{r,n}, ~~\gamma_{r,n,1,h}=\gamma_{r,n+1}=\gamma_{r,n+1,0,h}$ and $ \gamma_{r,n+1,1,h}=\gamma_{r,n+2},$ where $S_{n,h}$ is given by (\ref{1.213}) with $F_{r,n}=C_{r,n}^{T}$. \\
Note that when $h$ is small enough, $(I+shC_{r,n}^{T})^{-1}$  make sense,  and $S_{r,n,s,h}$ are symplectic matrices for $0\le s\le 1$ under the situation $B_{r,n}^T=B_{r,n}$ from Lemma 2.1.}

{\it {\bf Definition 2.6.} Let $\{\gamma_{n}\}$ possess the property (\ref{1.42}) with $\gamma_{0}=I_{2m}$, the connecting orbits $\gamma_{n,s,h}$ is given by (\ref{1.39}), where $n\in\bf N$. Denote by $\gamma_{d,h,n}$ a sequence of discrete symplectic matrices $\gamma_{0}, \gamma_{1}, \dots, \gamma_{n}$ and by $\gamma_{c,h,n}$ the joint path of a family of continuous symplectic matrices $\{\gamma_{0,s,h}, \gamma_{1,s,h}, \dots, \gamma_{(n-1),s,h}\}$ given by (\ref{1.39}), that is
\bea
\gamma_{c,h,n}=\gamma_{(n-1),\cdot,h}*\gamma_{(n-2),\cdot,h}*\dots*\gamma_{1,\cdot,h}*\gamma_{0,\cdot,h}.\lb{1.130}
\eea
Then $\gamma_{d,h,n}$ is called as the fundamental solution to the discrete Hamiltonian system (\ref{1.4}) with $B_{n}=B_{n}^{T}$, the maslov-type index $i_{\omega}(\gamma_{d,h,n})$ is defined as the Malov-type index $i_{\omega}(\gamma_{c,h,n})$ of $\gamma_{c,h,n}$, where $\omega \in {\bf U}$. }\\

{\it {\bf Remark 2.7. i).} The maslov-type index $i_{\omega}(\gamma_{d,h,n})$  is only dependent on $h, B_{l}$ of the discrete Hamiltonian system (\ref{1.4}) with $0 \le l \le n-1$, independent on $s$, hence it is well defined.\\
{\bf ii).} There is one to one correspondence between $\gamma_{d,h,n}$ and $\gamma_{c,h,n}$ from Lemma 2.1 and Lemma 2.2.\\
{\bf iii).} $\gamma_{d,h,N}$ can be seen as  $\tilde{\gamma}_{d,h,N} : \{0,\frac{1}{N},\dots,\frac{N-1}{N},\frac{N}{N}\} \rightarrow Sp(2m)$ with $\tilde{\gamma}_{d,h,N}(\frac{n}{N})=\gamma_{n}$ for $0 \le n \le N$. For notation simplicity, we still denote $\tilde{\gamma}_{d,h,N}$ by $\gamma_{d,h,N}$. Similarly, $\gamma_{c,h,N} : [0,1] \rightarrow Sp(2m)$ can be seen as $\gamma_{c,h,N}(t)=\gamma_{n,Nt-n,h}$ for $\frac{n}{N} \le t \le \frac{n+1}{N}$, where $\gamma_{n,Nt-n,h}$ is the same as $\gamma_{n,s,h}$ given by (\ref{1.39}) with $s$ replaced by $Nt-n$,  $0 \le n \le N-1$. Note that $\gamma_{c,h,N}(\frac{n}{N})=\gamma_{n}$ for $0\le n \le N.$
}

Based on Lemma 2.2, Definition 2.6 and Remark 2.7, we have\\

{\it {\bf Corollary 2.8.} For any continuous symplectic matrices path $\gamma(t)$ from the identity defined on $[0,1]$, there exist $h>0$ sufficiently small and the discrete Hamiltonian system (\ref{1.4}) with
\be
B_{n}=\left(
              \begin{array}{cc}
                A_{n} & C_{n} \\
                C_{n}^{T} & D_{n} \\
              \end{array}
            \right)\lb{1.286}
\ee
such that  $\gamma(t_{n})=\gamma_{d,h,N}(t_{n})=\gamma_{c,h,N}(t_{n})$ hold for $0 \le n \le N$, and  the corresponding  Maslov-type index and the nullity possess $i_{\omega}(\gamma)=i_{\omega}(\gamma_{d,h,N})=i_{\omega}(\gamma_{c,h,N}), \nu_{\omega}({\gamma})=\nu_{\omega}({\gamma_{d,h,N}})=\nu_{\omega}({\gamma_{c,h,N}})$, where $\gamma_{d,h,N}$ and $ \gamma_{c,h,N}$ are given by Definition 2.6 corresponding to the discrete Hamiltonian system with coefficient matrix (\ref{1.286}), $A_{n}^{T}=A_{n}, D_{n}^{T}=D_{n}$, $\omega \in {\bf U}, h=\frac{1}{N}, t_{n}=nh$.}

{\bf Proof.} We only need to take $\gamma_{d,h,N}(t_{n})$ as $\gamma(t_{n})$, thus we obtain the  existence of  $\gamma_{d,h,N}$, moreover we should have proved the Corollary 2.8 except $i_{\omega}(\gamma)=i_{\omega}(\gamma_{c,h,N})$ from Definition 2.6 if we had proved that $\gamma_{c,h,N}$ is well defined. The proof will be divided into two steps.

{\bf Step 1.} Well-definedness of $\gamma_{c,h,N}$.

Since the compactness of $t \in [0,1]$ and the continuity of $\gamma(t)$, then $\gamma(t)$ is uniformly continuous on $[0,1]$, thus for any $\epsilon >0$, there exists $h>0$ sufficiently small such that for $|t-t'|\le h$, we have
\bea
|\gamma(t)-\gamma(t')| <\epsilon.\lb{1.125}
\eea
In particular, take $t'=t_{n}$, it yields that for $|t-t_{n}|\le h$, there holds uniformly
\bea
|\gamma(t)-\gamma(t_{n})| <\epsilon.\lb{1.189}
\eea
Meanwhile, by the compactness of $t \in [0,1]$ and the continuity of $\gamma(t)$, then $\gamma(t)$ is bounded on $[0,1]$, together with the symplecticity of $\gamma(t)$,  hence there exist two positive constant $c_{1},c_{2}$ satisfying that
\bea
c_{1}<|\gamma(t)|<c_{2}, ~~~\text{for}~~t \in [0,1].\lb{1.190}
\eea
By (\ref{1.190}), it forces that
\bea
c_{1}<|\gamma(t_{n})|<c_{2}.\lb{1.191}
\eea
for $0 \le n \le N$.

By (\ref{1.42}), we have
\begin{eqnarray}
S_{n,h}=\gamma_{n+1}\gamma_{n}^{-1}:=\left(
               \begin{array}{cc}
               X_{n,11} & X_{n,12}\\
               X_{n,21} & X_{n,22}\\
               \end{array}
             \right).\lb{1.192}\end{eqnarray}
Hence
\bea
|S_{n,h}-I_{2m}|&=&|\gamma_{n+1}\gamma_{n}^{-1}-I_{2m}|\nn\\
                &\le&|\gamma_{n+1}-\gamma_{n}||\gamma_{n}^{-1}|\nn\\
                &\le&\frac{1}{c_{1}}\epsilon.\lb{1.193}
\eea
In the last line of (\ref{1.193}) we use (\ref{1.189}) and (\ref{1.191}).

Therefore, the existence of $B_{n}$ with components $A_{n}, C_{n}, D_{n}$ is obtained from Lemma 2.2. Moreover, by (\ref{1.204}), (\ref{1.205}), (\ref{1.206}) and (\ref{1.193}),  there holds
\bea
|hB_{n}|\le O(\epsilon), ~~~\text {for} ~~0 \le n \le N. \lb{1.197}
\eea
Then we have
\bea
|S_{n,s,h}-I_{2m}|\le O(\epsilon),\lb{1.198}
\eea
where $0 \le s \le 1 $, $S_{n,s,h}$ is given by (\ref{1.41}).

By Lemma 2.2, Remark 2.5, together with Definition 2.6, we obtain the corresponding $\gamma_{c,h,N}$ to $\gamma_{d,h,N}$ with $\gamma(t_{n})=\gamma_{d,h,N}(t_{n})=\gamma_{n}=\gamma_{c,h,N}(t_{n})$ for $0 \le n \le N$. Therefore $\gamma_{c,h,N}$ is well defined.

{\bf Step 2.} We need to verify $i_{\omega}(\gamma)=i_{\omega}(\gamma_{c,h,N})$.

For case $0\le t-t_{n}<h $, we have
\bea
|\gamma_{c,h,N}(t)-\gamma_{c,h,N}(t_{n})|&=& |\gamma_{n,Nt-n,h}-\gamma_{c,h,N}(t_{n})|\nn\\
                                         &=& |\gamma_{n,Nt-n,h}-\gamma_{n}|\nn\\
                                         &=& |S_{n,Nt-n,h}\gamma_{n}-\gamma_{n}|\nn\\
                                         &=& |S_{n,Nt-n,h}-I_{2m}||\gamma_{n}|\nn\\
                                         &\le& O(\epsilon).\lb{1.126}
\eea
We use iii) of Remark 2.7 in the second equality, the third equality uses (\ref{1.42}), the last inequality uses (\ref{1.191}) and (\ref{1.198}).

Similarly, for case $0\le t_{n}-t<h$, we also obtain
\bea
|\gamma_{c,h,N}(t)-\gamma_{c,h,N}(t_{n})|&=& |\gamma_{n-1,Nt-n+1,h}-\gamma_{c,h,N}(t_{n})|\nn\\
                                         &=& |\gamma_{n-1,Nt-n+1,h}-\gamma_{n}|\nn\\
                                         &=& |S_{n-1,Nt-n+1,h}\gamma_{n-1}-\gamma_{n-1,1,h}|\nn\\
                                         &=& |S_{n-1,Nt-n+1,h}\gamma_{n-1}-S_{n-1,1,h}\gamma_{n-1}|\nn\\
                                         &=& |S_{n-1,Nt-n+1,h}-S_{n-1,1,h}||\gamma_{n-1}|\nn\\
                                         &=& (|S_{n-1,Nt-n+1,h}-I_{2m}|+|S_{n-1,1,h}-I_{2m}|)|\gamma_{n-1}|\nn\\
                                         &\le& O(\epsilon).\lb{1.133}
\eea
We use iii) of Remark 2.7 in the second equality, the third equality holds due to (\ref{1.42}) and i) of  Remark 2.6, we use (\ref{1.42}) in the fourth equality, the last inequality uses (\ref{1.191}) and (\ref{1.198}).

On the other hand, for any $t \in [0,1]$, there exists $t_{n}$ satisfying $|t-t_{n}|<h$, hence there holds
\bea
|\gamma(t)-\gamma_{c,h,N}(t)| &=&|\gamma(t)-\gamma(t_{n})+\gamma(t_{n})-\gamma_{c,h,N}(t_{n})+\gamma_{c,h,N}(t_{n})-\gamma_{c,h,N}(t)|\nn\\
                              &\le&|\gamma(t)-\gamma(t_{n})|+|\gamma(t_{n})-\gamma_{c,h,N}(t_{n})|+|\gamma_{c,h,N}(t_{n})-\gamma_{c,h,N}(t)|\nn\\
                              &\le& O(\epsilon), \lb{1.127}
\eea
where the last inequality holds by (\ref{1.125}), (\ref{1.126}), (\ref{1.133}) and $\gamma(t_{n})=\gamma_{c,h,N}(t_{n})$. Together with Definition 2.6 and Remark 2.7,  the proof is complete.\hfill\hb\\

{\it {\bf Remark 2.9.} In general, $\gamma(t)$ doesn't wholly coincide with $\gamma_{c,h,N}(t)$ for $0 \le t \le 1$. }\\

Floquet theory is also valid for the discrete Hamiltonian system, we obtain\\

{\it {\bf Lemma 2.10.} Suppose that $\gamma_{d,h,N}$ be the fundamental solution to the discrete Hamiltonian system (\ref{1.4}) where $B_{n}$ satisfies the condition (\ref{1.291}), then $m^{0}_{h,\omega}=\nu_{\omega}(\gamma_{d,h,N}(1))$, where $\omega \in {\bf U}$.}

{\bf Proof.} Since the periodicity of the coefficient matrix $B_{n}$, there holds the periodicity of $S_{n,h}$ of the right hand side in (\ref{1.37}) with $F_{n}^T=C_{n}$, that is $S_{n,h}=S_{n+N,h}$, hence we have $\gamma_{d,h,N}(1+\frac{n}{N})=\gamma_{d,h,N}(\frac{n}{N})\gamma_{d,h,N}(1)$, so the result holds.\hfill\hb

\setcounter{equation}{0}
\section{The eigenvalues and eigenvectors of discrete Hamiltonian system}

In this section we discuss the eigenvalues and eigenvectors of the discrete Hamiltonian system under the boundary condition $z_{n+N}=\omega z_{n}$ for $\omega=e^{\sqrt {-1}\alpha}$ with $\alpha \in [0, 2\pi)$.

For the discrete Hamiltonian system (\ref{1.4}) with $h=1$,  we want to consider the following eigenvalue problem
\begin{eqnarray}
J(z_{n+1}-z_{n})=\lambda {\tilde z}_{n},\lb{1.23}
\end{eqnarray}
where the constant $\lambda \in {\bf C}$ will be determined.

The expression (\ref{1.23}) can be written as
\begin{eqnarray}
\left(
\begin{array}{cc}
y_{n}-y_{n+1}\\
x_{n+1}-x_{n}\\
\end{array}
\right)=
\lambda
\left(
\begin{array}{cc}
x_{n+1}\\
y_{n}\\
\end{array}
\right).\nn
\end{eqnarray}
That is
\bea
\left\{\begin{array}{ll}
               y_{n}-y_{n+1}=\lambda x_{n+1},~~~~~~~~~~~~~~~~~~~~~~~~~~~~~~~~~~~~ \\
               x_{n+1}-x_{n}=\lambda y_{n}. ~~~~~~~~~~~~~~~~~~~~~~~~~~~~~~~~~~~~~~ \end{array}\right. \nn\eea
Thus there holds
\bea
x_{n+2}=x_{n+1}+\lambda y_{n+1}=x_{n+1}+\lambda[y_{n}-\lambda x_{n+1}]=(1-{\lambda}^{2})x_{n+1}+\lambda y_{n},\nn
\eea
it yields that
\begin{eqnarray}
\left(
\begin{array}{cc}
x_{n+2}\\
y_{n+1}\\
\end{array}
\right)=
\left(
\begin{array}{cc}
(1-{\lambda}^{2})I_{m} & \lambda I_{m} \\
-\lambda I_{m} & I_{m}
\end{array}
\right)
\left(
\begin{array}{cc}
x_{n+1}\\
y_{n}\\
\end{array}
\right).\lb{1.24}
\end{eqnarray}
Let
\begin{eqnarray}
F(\lambda)=
\left(
\begin{array}{cc}
(1-{\lambda}^{2})I_{m} & \lambda I_{m} \\
-\lambda I_{m} & I_{m}\\
\end{array}
\right),\lb{1.25}
\end{eqnarray}
we want to find $f,\lambda$ and ${\mathcal{G}}=({\mathcal{G}}_{1}^{T},{\mathcal{G}}_{2}^{T})^{T}$ with ${\mathcal{G}}_{1},{\mathcal{G}}_{2} \in  {\bf C}^{m}$ satisfying $F(\lambda){\mathcal{G}}=f{\mathcal{G}}$, where ${\mathcal{G}} \neq 0$.

Firstly, it is deduced that ${\tilde z}_{n}=f^{n}{\mathcal{G}}$, ${\tilde z}_{n+1}=f^{n+1}{\mathcal{G}}$ is a solution to (\ref{1.24}) from the fact that $F(\lambda){\mathcal{G}}=f{\mathcal{G}}$, moreover it is the eigenvector of (\ref{1.23}), thus $f^{N}=\omega$ from $z_{n+N}=\omega z_{n}$. To be precise, $f_{k}=e^{\frac{(2k\pi+\alpha) \sqrt{-1}}{N}}$ with $k=0,1,\dots,N-1$. Write $f_{k}=e^{\sqrt{-1}{2}\alpha _{k}}$, where $\alpha _{k}= \frac{k\pi+\frac{\alpha}{2}}{N}$.

Secondly, based on the fact that $F(\lambda){\mathcal{G}}=f{\mathcal{G}}$ with ${\mathcal{G}} \neq 0$, there holds $|F(\lambda)-f_{k} I_{2m}|=0$, together with (\ref{1.25}), we have $(1-{\lambda}^2 -f_{k})(1-f_{k})+{\lambda}^2=0$, thus ${\lambda}^2=-\frac{1-2f_{k} +f_{k}^{2}}{f}$, combining with $f_{k}=e^{\sqrt{-1}{2}\alpha _{k}}$, we get that $\lambda_{k}^{\pm} =\pm \sqrt {-1}\frac{1-f_{k}}{\sqrt{f_{k}}}=\pm \sqrt {-1}\frac{1-e^{\sqrt{-1}{2}\alpha _{k}}}{e^{\sqrt{-1}\alpha _{k}}}=\pm 2sin\alpha_{k}$ with $m$ multiplicity respectively. Substituting it into $F(\lambda){\mathcal{G}}=f{\mathcal{G}}$, where ${\mathcal{G}}=({\mathcal{G}}_{1}^{T},{\mathcal{G}}_{2}^{T})^{T}$, we obtain the following results:

When $\lambda_{k}=2sin\alpha_{k}$ with $k \in \{0,\dots,N-1\}$. By simple computations, we have ${\mathcal{G}}=(-\sqrt{-1}e^{\sqrt{-1}\alpha_{k}}a_{0}^{T},a_{0}^{T})^{T}$, where $a_{0}\in {\bf R}^{m}$.  Hence
\begin{eqnarray}
{\tilde z}_{k,n}&=&f_{k}^{n}{\mathcal{G}}=(e^{\sqrt{-1}2\alpha_{k}})^{n}(-\sqrt{-1}e^{\sqrt{-1}\alpha_{k}}a_{0}^{T},a_{0}^{T})^{T}\nn\\
&=&([sin(2n+1)\alpha_{k}-\sqrt{-1}cos(2n+1)\alpha_{k}]a_{0}^{T},[cos2n\alpha_{k}+\sqrt{-1}sin2n\alpha_{k}]a_{0}^{T})^{T}.\lb{1.68}
\end{eqnarray}

When $\lambda_{k}=-2sin\alpha_{k}$ with $k \in \{0,\dots,N-1\}$. We have ${\mathcal{G}}=(\sqrt{-1}e^{\sqrt{-1}\alpha_{k}}a_{0}^{T},a_{0}^{T})^{T}$, where $a_{0}\in {\bf R}^{m}$.  Hence
\begin{eqnarray}
{\tilde z}_{k,n}&=&f_{k}^{n}{\mathcal{G}}=(e^{\sqrt{-1}2\alpha_{k}})^{n}(\sqrt{-1}e^{\sqrt{-1}\alpha_{k}}a_{0}^{T},a_{0}^{T})^{T}\nn\\
&=&([-sin(2n+1)\alpha_{k}+\sqrt{-1}cos(2n+1)\alpha_{k}]a_{0}^{T},[cos2n\alpha_{k}+\sqrt{-1}sin2n\alpha_{k}]a_{0}^{T})^{T}.\lb{1.69}
\end{eqnarray}

In particular, we further consider the periodic boundary case $\omega =1$, the case $\omega =-1$ is similarly dealt with. In this case, we write $\alpha_{k}$ as $\theta_{k}$.

{\bf Case 1.} $f_{k}=e^{\sqrt{-1}2\theta_{k}}, \lambda_{k}=2sin\theta_{k}$ with $\theta _{k}= \frac{k\pi}{N}$, where $k \in \{0,\dots,N-1\}$. By simple computations, we have ${\mathcal{G}}=(-\sqrt{-1}e^{\sqrt{-1}\theta_{k}}a_{0}^{T},a_{0}^{T})^{T}$, where $a_{0}\in {\bf R}^{m}$.  Hence
\begin{eqnarray}
{\tilde z}_{k,n}&=&f_{k}^{n}{\mathcal{G}}=(e^{\sqrt{-1}2\theta_{k}})^{n}(-\sqrt{-1}e^{\sqrt{-1}\theta_{k}}a_{0}^{T},a_{0}^{T})^{T}\nn\\
&=&([sin(2n+1)\theta_{k}-\sqrt{-1}cos(2n+1)\theta_{k}]a_{0}^{T},[cos2n\theta_{k}+\sqrt{-1}sin2n\theta_{k}]a_{0}^{T})^{T}.\lb{1.26}
\end{eqnarray}
Since $\theta_{k}+\theta_{N-k}=\pi$, then $sin\theta_{k}=sin\theta_{N-k}$, thus $\lambda_{k}=\lambda_{N-k}$. So this case falls into the following three subcases.\\
{\bf Case 1.1.} $k \in \{1,2,\dots, [\frac{N}{2}]\}$. \\
$\tilde{z}_{k,n}$ has the form (\ref{1.26}), and note that when $N$ is even, $\theta_{\frac{N}{2}}=\frac{\pi}{2}$, $\lambda_{\frac{N}{2}}=2$, and
\begin{eqnarray}
{\tilde z}_{\frac{N}{2},n}=(-1)^{n}(a_{0}^{T},a_{0}^{T})^{T}.\lb{1.27}
\end{eqnarray}
{\bf Case 1.2.} $k \in \{[\frac{N}{2}]+1,\dots, N-1\}$. \\
Since $\theta_{k}+\theta_{N-k}=\pi$, thus $\lambda_{k}=\lambda_{N-k}$. Moreover,
\begin{eqnarray}
{\tilde z}_{k,n}&=&([sin(2n+1)\theta_{k}-\sqrt{-1}cos(2n+1)\theta_{k}]b_{0}^{T},[cos2n\theta_{k}+\sqrt{-1}sin2n\theta_{k}]b_{0}^{T})^{T}\nn\\
              &=&([sin(2n+1)\theta_{N-k}+\sqrt{-1}cos(2n+1)\theta_{N-k}]b_{0}^{T},[cos2n\theta_{N-k}-\sqrt{-1}sin2n\theta_{N-k}]b_{0}^{T})^{T}.\lb{1.28}
\end{eqnarray}
Based on the discussions above, the eigenvector $\tilde{z}_{k,n}$ associated with $\lambda_{k}=2sin\theta_{k}$ can be written as
\begin{eqnarray}
\tilde z_{k,n}=
\left(
\begin{array}{cc}
sin(2n+1)\theta_{k}\hspace*{1.5mm}a_{k}\\
cos2n\theta_{k}\hspace*{1.5mm}a_{k}
\end{array}
\right)\lb{1.29}
\end{eqnarray}
and
\begin{eqnarray}
\tilde z_{k,n}=
\left(
\begin{array}{cc}
-cos(2n+1)\theta_{k}\hspace*{1.5mm}b_{k}\\
sin2n\theta_{k}\hspace*{1.5mm}b_{k}
\end{array}
\right),\lb{1.30}
\end{eqnarray}
where $a_{k}, b_{k} \in {\bf R}^{m}$ and $k \in \{1,2,\dots,[\frac{N}{2}]\}$.\\
{\bf Case 1.3.} $k=0$.\\
In this case, $\lambda_{0}=0$, ${\mathcal{G}}=({\mathcal{G}}_{1}^{T},{\mathcal{G}}_{2}^{T})^{T}$, where ${\mathcal{G}}_{1},{\mathcal{G}}_{2} \in {\bf R}^{m}$. Moreover $\tilde{z}_{0,n}={\mathcal{G}}$.

Similar arguments yield that

{\bf Case 2.} $f_{k}=e^{\sqrt{-1}2\theta_{k}}, \lambda_{k}=-2sin\theta_{k}$ with $\theta _{k}= \frac{k\pi}{N}$, where $k \in \{0,\dots,N-1\}$. Similar to Case 1, the eigenvector $\tilde{z}_{k,n}$ can be written as
\begin{eqnarray}
\tilde z_{k,n}=
\left(
\begin{array}{cc}
-sin(2n+1)\theta_{k}\hspace*{1.5mm}c_{k}\\
cos2n\theta_{k}\hspace*{1.5mm}c_{k}
\end{array}
\right)\lb{1.31}
\end{eqnarray}
or
\begin{eqnarray}
\tilde z_{k,n}=
\left(
\begin{array}{cc}
cos(2n+1)\theta_{k} \hspace*{1.5mm}d_{k}\\
sin2n\theta_{k} \hspace*{1.5mm}d_{k}
\end{array}
\right).\lb{1.32}
\end{eqnarray}
In particular, when $N$ is even, then $\theta_{\frac{N}{2}}=\frac{\pi}{2}$, $\lambda_{\frac{N}{2}}=-2$, and
\begin{eqnarray}
{\tilde z}_{\frac{N}{2},n}=(-1)^{n}(-a_{0}^{T},a_{0}^{T})^{T},\lb{1.36}
\end{eqnarray}
where $c_{k},d_{k}\in {\bf R}^{m}$ and $k \in \{1,2,\dots,[\frac{N}{2}]\}$, $a_{0} \in {\bf R}^{m}$.

\setcounter{equation}{0}
\section{ Morse index and Maslov-type index}
As usual, the symplectic group $Sp(2m)$ is defined by\[Sp(2m) = \{M
\in GL(2m,{\bf R}) \mid M^{T}JM = J\},\] whose topology is induced
from that of ${\bf R}^{4m^2}$. For $\tau>0$ we are interested in
paths in $Sp(2m)$:
\begin{eqnarray}
&&\mathcal {P}_{\tau}(2m) = \{\gamma \in C([0,
\tau], Sp(2m)) \mid \gamma(0) = I_{2m}\},\nn\\
&&{\hat {\mathcal {P}}}_{\tau}(2m) = \{\gamma \in C^{1}([0,
\tau], Sp(2m)) \mid \gamma(0) = I_{2m}, {\dot{\gamma}}(\tau)={\dot{\gamma}}(0)\gamma(\tau)\},\nn
\end{eqnarray}
which are equipped with
the topology induced from that of $Sp(2m)$. The following real
function was introduced in \cite{Lon4}:\[D_{\omega}(M) =
(-1)^{n-1}\bar{\omega}^n det(M-\omega I_{2m}),~\forall~\omega \in {\bf
U}, M \in Sp(2m).\]where ${\bf
U}$ is the unit circle in the complex plane. Thus for any $\omega\in {\bf U}$ the following
codimension $1$ hypersurface in $Sp(2m)$ is defined in
\cite{Lon4}:\[Sp(2m)^0_{\omega} = \{M \in Sp(2m) \mid
D_{\omega}(M) = 0\}.\] For any $M \in Sp(2m)^0_{\omega}$, we
define a co-orientation of $Sp(2n)^0_{\omega}$ at $M$ by the
positive direction $\frac{d}{dt}M e^{t\epsilon J} |_{t=0}$ of the
path $M e^{t\epsilon J}$ with $0\leq t \leq 1$ and $\epsilon > 0$
being sufficiently small. Let
\begin{eqnarray}
&&Sp(2n)_{\omega}^* = Sp(2m)\setminus Sp(2m)^0_{\omega},\nn\\
&&\mathcal{P}^*_{\tau,\omega}(2m) = \{\gamma \in \mathcal {P}_{\tau} (2m)\mid \gamma_{\tau}\in Sp(2m)^*_{\omega}\},\nn\\
&&{\hat {\mathcal{P}}}^*_{\tau,\omega}(2m) = \{\gamma \in {\hat {\mathcal{P}}}_{\tau} (2m)\mid \gamma_{\tau}\in Sp(2m)^*_{\omega}\},\nn\\
&&\mathcal {P}^0_{\tau,\omega}(2m)=\mathcal {P}_{\tau} (2m)\setminus \mathcal{P}^*_{\tau,\omega}(2m),\nn\\
&&{\hat {\mathcal {P}}}^0_{\tau,\omega}(2m)={\hat {\mathcal{P}}}_{\tau} (2m)\setminus {\hat {\mathcal {P}}}^*_{\tau,\omega}(2m).\nn\end{eqnarray}
For any two continuous arcs $\xi$ and $\eta : [0, \tau ] \longrightarrow Sp(2m)$ with $\xi(\tau) = \eta(0)$, it is defined as usual:\\
$${\eta}*{\xi}(t)=\left\{ \begin{array}{ll}\xi(2t),~if~0\leq t\leq \tau/2,\\
\eta(2t-1),~if~\tau/2\leq t \leq \tau.\end{array}\right.$$
Given any two $2m_k \times 2m_k$ matrices of square block form $M_{k}=\left(
              \begin{array}{cc}
                A_k & B_k \\
                C_k & D_k\\
              \end{array}
            \right)$
with $k = 1, 2$, as in \cite{Lon4}, the $\diamond$-product of $M_1$
and $M_2$ is defined by the following $2(m_1 +m_2) \times 2(m_1
+m_2)$ matrix $M_1 \diamond M_2$:\\$$M_1 \diamond M_2=\left(
                            \begin{array}{cccc}
                              A_1& 0 & B_1 & 0 \\
                              0 & A_2 & 0 & B_2 \\
                              C_1 & 0 & D_1 & 0 \\
                              0 & C_2 & 0 & D_2 \\
                            \end{array}
                          \right)$$
Denote by $M^{\diamond k}$ the $k$-fold $\diamond$-product $M
\diamond \cdots \diamond M$. Note that the $\diamond$-product of any
two symplectic matrices is symplectic. For any two paths $\gamma_{j}
\in \mathcal {P}_{\tau}(2n_j)$ with $j = 0$ and $1$, let $\gamma_1
\diamond \gamma_2(t)=\gamma_1(t)\diamond \gamma_2(t)$ for all $t \in
[0, \tau]$.

{\it {\bf Definition 4.1}(cf. \cite{Lon4}) For any $\omega \in {\bf U}$
and $M \in Sp(2m)$, we define
\begin{eqnarray}\nu_{\omega}(M) = dim_{{\bf C}} ker_{{\bf
C}}(M - \omega I_{2m}).
\end{eqnarray}}
For any $\tau > 0$ and $\gamma \in \mathcal {P}_{\tau} (2m)$,
define
\begin{eqnarray}
\nu_{\omega}(\gamma) = \nu_{\omega}(\gamma(\tau )).\end{eqnarray}
If $\gamma \in
\mathcal {P}^*_{\tau,\omega}(2m)$,
define
\begin{eqnarray}i_{\omega}(\gamma) = [Sp(2m)^0_{\omega}:
\gamma\ast \xi_m],\lb{1.295}\end{eqnarray}
where the right hand side of
$(\ref{1.295})$ is the usual homotopy intersection number, and the
orientation of $\gamma\ast \xi_m$ is its positive time direction
under homotopy with fixed end points.

If $\gamma\in \mathcal {P}^0_{\tau,\omega}(2m)$, we let $\mathcal
{F}(\gamma)$ be the set of all open neighborhoods of $\gamma$ in
$\mathcal {P}_{\tau}(2m)$, and define
\begin{eqnarray}
i_{\omega}(\gamma) = \sup_{ U\in \mathcal {F}(\gamma)}{
\inf{\{i_{\omega}(\beta) \mid \beta \in U\cap\mathcal
{P}^*_{\tau,\omega}(2m)\}}}.
\end{eqnarray}
Then\begin{eqnarray}(i_{\omega}(\gamma), \nu_{\omega}(\gamma))
\in {\bf Z} \times \{0, 1, \cdots , 2m\} \end{eqnarray} is called
the index function of $\gamma$ at $\omega$.

Let $\Omega^0(M)$ be the path connected component containing $M =
\gamma(\tau)$ of the set\begin{eqnarray} \Omega(M) = \{N \in Sp(2n)
\mid  \sigma(N)
\cap {\bf U} = \sigma(M) \cap {\bf U}~and~~~~\nn\\
\nu_{\lambda}(N)=\nu_{\lambda}(M),\forall \lambda \in \sigma(M) \cap
{\bf U}\}\nn\end{eqnarray} Here $\Omega^0(M)$ is called the {\it
homotopy component} of $M$ in $Sp(2n)$.

In \cite{Lon4}, the following symplectic matrices were introduced as
basic normal forms: \begin{eqnarray}D(\lambda)&=&\left(
                                                    \begin{array}{cc}
                                                      \lambda & 0 \\
                                                      0 & \lambda^{-1} \\
                                                    \end{array}
                                                  \right),
                                                  \lambda=\pm 2,\nn\\
N_1(\lambda,b)&=&\left(
 \begin{array}{cc}
  \lambda & b \\
   0& \lambda \\
   \end{array}
   \right),\lambda=\pm1,b=\pm1,0,\nn\\
R(\theta)&=&\left(
            \begin{array}{cc}
              \cos{\theta} & -\sin{\theta} \\
              \sin{\theta} & \cos{\theta} \\
            \end{array}
          \right),\theta\in (0,\pi)\cup(\pi,2\pi),\nn\\
N_2(\omega,B)&=&\left(
                \begin{array}{cc}
                  R(\theta) & B \\
                  0 &  R(\theta) \\
                \end{array}
              \right),\theta\in (0,\pi)\cup(\pi,2\pi),
\nn\end{eqnarray}
where $B=\left(
           \begin{array}{cc}
             b_1 & b_2 \\
             b_3 & b_4 \\
           \end{array}
         \right)
$ with $b_i\in {\bf R}$ and $b_2\neq b_3$.

Next we define the standard smooth paths ${\hat{\beta}}_{j}$ in $\mathcal {P}_{1}(2m)$ and the corresponding coefficient matrices ${\hat{B}}_{j}$ which were introduced in (\cite{LoZ1}) according to case $m=1$ or $m\geq 2$ as follows (or cf. section 5.3 of \cite{Lon4} ).

Case $m=1$. \\
Set $\omega(t)=\frac{1}{2}(1+cos(2\pi t))$,

\be \omega_{1}(t) = \left\{\begin{array}{ll}
                    \omega(t), & {\rm if\;}0\le t\le \frac{1}{2},  \\
                      0, & {\rm \frac{1}{2}\le t\le 1}.\\
                      \end{array}
\right.\lb{1.100}  \ee

\be \omega_{1}(t) = \left\{\begin{array}{ll}
                       0, & {\rm if\;}0\le t\le \frac{1}{2},  \\
                      \omega(t), & {\rm \frac{1}{2}\le t\le 1}. \\
                      \end{array}
\right.\lb{1.101}  \ee

and
\begin{eqnarray}K=\left(
              \begin{array}{cc}
                0 & 1 \\
                1 & 0 \\
              \end{array}
            \right).\lb{1.109}
\end{eqnarray}
Define
\be
{\hat{\beta}_{0}}(t)=D(2^{t}),~~ {\hat{B}_{0}}=(ln2)K.\lb{1.102}
\ee
If $j\in 2{\bf Z}+1$, we define
\be
{\hat{\beta}_{j}}(t)=R(j\pi t),~~ {\hat{B}_{j}}=j\pi I_{2m}.\lb{1.103}
\ee
If $j\in 2{\bf Z}\setminus \{0\}$, we define
\be
{\hat{\beta}_{j}}(t)=D(2^{\omega_{2}(t)})R(1-\omega_{1}(t)j\pi),
~~ {\hat{B}_{j}}(t)=-j\pi{\dot {\omega}_{1}(t)}I_{2m} +(ln2){\dot {\omega}_{2}(t)}K.\lb{1.104}
\ee
Case $m\ge 2$.\\
Define $2m\times 2m$ matrices
\begin{eqnarray}
X_{j}&=&diag(0,(j-m+2)\pi, \pi,\cdots, \pi),\nn\\
Y    &=&diag(ln2,0,\cdots, 0),\nn\\
Z_{j}&=&diag((j-m+2)\pi, \pi,\cdots, \pi).\nn
\end{eqnarray}
If $(-1)^{m+j}=-1$, we define
\begin{eqnarray}{\hat {\beta}_{j}}(t)&=&exp(tJ{\hat {B}_{j}}),\lb{1.105}\\
{\hat {B}_{j}}&=&\left(
 \begin{array}{cc}
  X_{j} & Y \\
  Y     & X_{j} \\
   \end{array}
   \right).\lb{1.106}
\end{eqnarray}
If $(-1)^{m+j}=1$, we define
\begin{eqnarray}{\hat {\beta}_{j}}(t)&=&exp(tJ{\hat {B}_{j}}),\lb{1.107}\\
{\hat {B}_{j}}&=&\left(
 \begin{array}{cc}
  Z_{j} & 0 \\
  0     & Z_{j} \\
   \end{array}
   \right).\lb{1.108}
\end{eqnarray}
There are other choices of ${\hat{\beta}}_{j}$ and ${\hat{B}}_{j}$ which will appear in the Case 2 and Case 3 of Section 6, we omit them here.

As (5.4.6) in \cite{Lon4}, we define a small perturbation of $\gamma_{c,h,N}$
\be
\gamma_{c,s}(t)=\gamma_{c,h,N}(t)P^{-1}R_{m_{1}}(s\rho(t))\theta_{0})\dots R_{m_{p+2q}}(s\rho(t)\theta_{0})P, \lb{1.44}
\ee
satisfying $\gamma_{c,s}(t)$ converges to $\gamma_{c,h,N}(t)$ in $C^{1}([0,\tau],Sp(2m))$ as $s \rightarrow 0$ and $\gamma_{c,0}(t)=\gamma_{c,h,N}(t)$, where $-1 \le s \le 1$,  $P$ is some symplectic matrix, $\rho \in C^{2}([0,\tau],[0,1])$ with some properties.

To simplify the notation, we write the perturbation of $\gamma_{d,h,N}$ as $\gamma_{d,s}$. Note that for any $0\le n \le N$, there holds $\gamma_{c,h,N}(\frac{n}{N})=\gamma_{d,h,N}(\frac{n}{N})$, then $\gamma_{c,0}(\frac{n}{N})=\gamma_{d,h,N}(\frac{n}{N})$, thus $\gamma_{c,s}(\frac{n}{N})$ is a perturbation of $\gamma_{d,h,N}(\frac{n}{N})$, that is
\be
\gamma_{d,s}(\frac{n}{N})=\gamma_{d,h,N}(\frac{n}{N})P^{-1}R_{m_{1}}(s\rho(\frac{n}{N}))\theta_{0})\dots R_{m_{p+2q}}(s\rho(\frac{n}{N})\theta_{0})P. \lb{1.45}
\ee
Note that $\gamma_{d,s}(1) \in Sp(2n)_{\omega}^*$ for $s\in [-1,1]\setminus\{0\}$, and we have the following property about $\gamma_{d,s}$\\

{\it {\bf Lemma 4.2.} Suppose that $\gamma_{d,h,N}$ be a sequence of symplectic matrices, and $\gamma_{d,h,N}(1) \in Sp^{0}_{\omega}(2m)$ for $\omega \in {\bf U},$  $\gamma_{d,s}(\cdot)$ is given by (\ref{1.45}), there holds
\bea
i_{\omega}(\gamma_{d,s})-i_{\omega}(\gamma_{d,-s})=\nu_{\omega}(\gamma_{d,h,N}),
\eea
where $ 0<s\le 1$.}

{\bf Proof.} Since $i_{\omega}(\gamma_{d,s})=i_{\omega}(\gamma_{c,s}), i_{\omega}(\gamma_{d,-s})=i_{\omega}(\gamma_{c,-s}), \nu_{\omega}(\gamma_{d,h,N})=\nu_{\omega}(\gamma_{c,h,N})$, together with Lemma 2.2 and
\bea
i_{\omega}(\gamma_{c,s})-i_{\omega}(\gamma_{c,-s})=\nu_{\omega}(\gamma_{c,h,N}).
\eea
cf. Theorem 5.4.1 of \cite{Lon4}, the result is followed.\hfill\hb

Next we introduce the following two bilinear form on $\tilde W$:

The first one is $\mathcal{A}_{h,\omega}(\tilde {z},\tilde {w})$ given by (\ref{1.47}), where $B_{n}$ satisfies the condition (\ref{1.291}). Denote by $m_{h,\omega}^{-}, m_{h,\omega}^{0}, m_{h,\omega}^{+}$ the negative, null, positive Morse indices of $\mathcal{A}_{h,\omega}$.

The other one is
\bea
\mathcal{A}_{r,h,\omega}(\tilde {z},\tilde {w})=\Sigma_{n=1}^{N}({-J\frac{{z}_{n+1}-z_{n}}{h},\tilde{w}_{n}})-(B_{r,n}\tilde{z}_{n}, \tilde{w}_{n}).\lb{1.48}
\eea
 also, denote by $m_{r,h,\omega}^{-}, m_{r,h,\omega}^{0}, m_{r,h,\omega}^{+}$ the dimension of the negative, null, positive invariant subspace of $\mathcal{A}_{r,h,\omega}$ respectively, where $B_{r,n}$ are real symmetric $2m\times 2m$ matrices with $B_{r,n+N}=B_{r,n}, ~~0 \le r \le 1$. Moreover we can introduce\\

{\it{\bf Definition 4.3.} The splitting numbers $\mathcal{S}_{h,\omega}^{\pm}$ of the discrete Hamiltonian system (\ref{1.4}) with $B_{n}=B_{n}^{T}$ are defined by
\bea
\mathcal{S}_{h,\omega}^{\pm}=\lim _{\theta \rightarrow 0^{\pm}}m_{h,{\pm\theta}\omega e^{\sqrt -1}}^{-}-m_{h,\omega}^{-}.\lb{1.61}
\eea}\\

{\it {\bf Lemma 4.4.} Suppose that $B_{r,n}=B(r,\frac{n}{N})$ is continuous with respect to $r$, where $B:[0,1]\times \{0,\frac{1}{N},\dots,\frac{N}{N}\} \rightarrow {\mathcal L}_{s}({\bf R}^{2m})$ with $B_{0,n}=B_{n}$, $0 \le n \le N$. Provided that $m_{r,h,\omega}^{0}=0, \forall r \in [0,1]$, there hold $$m_{r,h,\omega}^{-}=m_{0,h,\omega}^{-}=m_{h,\omega}^{-}, ~~m_{r,h,\omega}^{+}=m_{0,h,\omega}^{+}=m_{h,\omega}^{+}, ~~\forall r \in [0,1].$$}
{\bf Proof.} Choose two close enough numbers $r,r^{'} \in [0,1]$, there hold $m_{r,h,\omega}^{-} \le m_{r^{'},h,\omega}^{-},~~~ m\xi(r,t_{n})_{r,h,\omega}^{+} \le m_{r^{'},h,\omega}^{+}$, hence we have
$2Nm=m_{r,h,\omega}^{-}+m_{r,h,\omega}^{0}+m_{r,h,\omega}^{+} \le m_{r^{'},h,\omega}^{-}+m_{r^{'},h,\omega}^{0}+m_{r^{'},h,\omega}^{+}=2Nm.$ Therefore $m_{r,h,\omega}^{-},m_{r,h,\omega}^{+}$ are locally constants, together with the connectedness of $[0,1]$, they are also globally constants.\hfill\hb\\

{\it{\bf Definition 4.5.} (cf. Definition 9.1.4 of \cite{Lon4}) 
For any $M \in Sp(2n)$ and $\omega \in {\bf U}$, {\it the splitting
numbers} $S^{\pm}_M (\omega)$ of $M$ at $\omega$ are defined
by
\bea
S_{M}^{\pm}({\omega})=\lim _{\theta \rightarrow 0^{\pm}}i_{{\pm\theta}\omega e^{\sqrt -1}}(\gamma)-i_{\omega}(\gamma),\lb{1.297}
\eea
for any path $\gamma \in \mathcal {P}_{\tau} (2n)$ satisfying
$\gamma(\tau) = M$.\\

{\it {\bf Lemma 4.6.}(cf. Lemma 9.1.5 of \cite{Lon4}) The splitting numbers $S_{M}^{\pm}({\omega})$ are independent of the path $\gamma \in \mathcal{P}_{\tau}(2m)$ appearing in (\ref{1.297}).}

\setcounter{equation}{0}
\section{Construction of the homotopy of discrete symplectic paths}
In this section, we suppose that $i_{\omega}(\gamma_{d,h,N})=j$, $\nu_{\omega}(\gamma_{d,h,N}(t_{N}))=0,$ where $j \in {\bf Z}, t_{n}=nh$ with $0 \le n \le N$, $h$ is some sufficiently small positive number.\\

{\it {\bf Definition 5.1} Denote by $\zeta(s,h,N)$ the sequence of $\{\zeta(s,h,N)(t_{n})\}$, where $0\le n \le N$.}

The motivation of the construction of the homotopy of discrete symplectic paths is that we want to reduce the proof of Theorem 1.1 for general discrete Hamiltonian system to the proof for the discrete Hamiltonian with coefficient matrix ${\hat B}_{j,n}$ which corresponds to ${\hat B}_{j}(t)$ associated to the standard path ${\hat \beta}_{j}(t)$, where ${\hat B}_{j,n}={\hat B}_{j}(t_{n})$. In this case the Morse indices of the corresponding functional $\mathcal{A}_{h,\omega}$ at the critical point $\tilde z=0$ can be easily computed.

Firstly we outline the construction of the homotopy of discrete symplectic paths.  The basic strategy is that by virtue of the homotopy $\xi(r,t)$ for the continuous case between ${\gamma}_{c,h,N}(t)$ and ${\hat \beta}_{j}(t)$ with $0 \le r, t \le 1$, we can obtain the corresponding homotopy $\zeta(r,h,N)(t_{n})$ for the discrete case between ${\gamma}_{c,h,N}(t_{n})$ and ${\hat \beta}_{j}(t_{n})$, that is the homotopy connecting ${\gamma}_{d,h,N}(t_{n})$ to ${\hat \beta}_{j}(t_{n})$. To reduce the proof of Theorem 1.1 to the proof for the discrete Hamiltonian with coefficient matrix ${\hat B}_{j,n}$ which corresponds to ${\hat B}_{j}(t)$ associated to the standard path ${\hat \beta}_{j}(t)$, we must guarantee that $\nu_{\omega}(\zeta(r,h,N)(t_{N}))=0$ for $0 \le r \le 1$.

Moreover note that ${\gamma}_{c,h,N}(t)$ is only continuous in $t$, while ${\hat {\beta}}_{j}(t)$ is $C^{1}$ in $t$, thus $\xi(r,t)$ can be divided into two parts. The first part $\xi(r,t)$  connecting ${\gamma}_{c,h,N}(t)$ to certain path $\psi \in {\hat {\mathcal {P}}}_{1}(2m)$ is continuous in $t$ for $0 \le r \le r_{1}$, the second part $\xi(r,t)$ between $\psi$ and ${\hat {\beta}}_{j}(t)$ is $C^{1}$ in $t$ for $r_{1} \le r \le 1$, where $r_{1} \in [0,1]$. For $\xi(r,t)$ with $0 \le r \le r_{1}$, since $\xi(r,t)$ is just continuous in $t$, we can apply Corollary 2.8 to discretize $\xi(r,t)$, furthermore yield the corresponding coefficient matrices $B_{r,n}$ of discrete Hamiltonian system (\ref{1.211}) from Lemma 2.2, hence finally obtain $\zeta(r,t_{n})$ as the fundamental solutions to the discrete Hamiltonian system (\ref{1.211}) with the coefficient matrices $B_{r,n}$, which connects ${\gamma}_{c,h,N}(t_{n})$ to $\psi(t_{n})$ for $0 \le n \le N$.  Note that $\zeta(r,h,N)(t_{n})=\xi(r,t_{n})$ for $0 \le r \le r_{1}.$

For the second part $\xi(r,t)$ with $r_{1} \le r \le 1$, in order to guarantee that the discrete Hamiltonian system with coefficient matrix ${\hat B}_{j,n}$, we cannot discretize $\xi(r,t)$ as the first part above. Fortunately, since $\xi(r,t)$ is $C^{1}$ in $t$, we can obtain the coefficient matrix $B_{r}(t)$ given by (\ref{1.110}) of continuous Hamiltonian system, moreover take $B_{r,n}=B_{r}(t_{n})$ as the coefficient matrices of discrete Hamiltonian system, finally $\zeta(r,h,N)(t_{n})$ is obtained as the fundamental solutions to discrete Hamiltonian system with coefficient matrices $B_{r,n}$. In general, $\zeta(r,h,N)(t_{n})\ne$ for $r_{1} \le r \le 1.$ In particular, $\psi$ (i.e.$\xi(r_{1},t)$) has been discretized twice, that is in the first part and the second part respectively, thus there is a gap between two different discretization paths, finally we obtain the third homotopy connecting two different discretization paths. Denote by ${\hat{\gamma}}_{j,h,d,N}$ the fundamental solution to the discrete Hamiltonian system with coefficient matrix ${\hat B}_{j,n}$.

The construction is divided into four steps. Step 1 briefly sketches the homotopy for the continuous case, Step 2-Step 4 is the homotopy for the discrete case. To be precise, Step 2 is devoted to the discrete homotopy corresponding to $\xi(r,t)$ with $0 \le r \le r_{1}$, Step 3 is used to the construction of the homotopy of the gap between two different discretizations of $\xi(r_{1},t)$, Step 4 is contributed to the discrete homotopy corresponding to $\xi(r,t)$ for $r_{1} \le r \le 1$.

{\bf Step 1.} The homotopy for the continuous case

Since $\gamma_{c,h,N} \in \mathcal {P}_{1}(2m)$, by the method of Step 2 in the proof of Theorem 6.1.8 in \cite{Lon4}, we can
choose a path $\psi \in {\hat {\mathcal {P}}}_{1}(2m)$ such that $\gamma_{c,h,N}\sim \psi$ with fixed end points. By Lemma 6.1.5 in \cite{Lon4}, there exists a homotopy between $\psi$ and ${\hat{\beta}}_{j}$ in ${\hat {\mathcal {P}}}_{1,\omega}(2m)$. Summing up, together with $\nu_{\omega}(\gamma_{c,h,N}(1))=\nu_{\omega}(\gamma_{d,h,N}(1))=0$, there exists a homotopy $\xi(r,t):[0,1]^{2} \rightarrow Sp(2m)$ connecting ${\gamma}_{c,h,N}(t)$ to $\psi$ for $0\le r \le r_{1}$ with
\bea
\xi(0,t)=\gamma_{c,h,N}, ~~~\xi(r_{1},t)=\psi \lb{1.288}
\eea
and connecting $\psi$ to certain standard path ${\hat{\beta}}_{j}$ in $\hat{\mathcal P}_{1,\omega}^{*}(2m)$ for $r_{1} \le r \le 1$ with
\bea
\xi(r_{1},t)=\psi, ~~~\xi(1,t)={\hat{\beta}}_{j},\lb{1.289}
 \eea
 which is continuous in  $r,t  \in [0,1]$, where $r_{1} \in (0,1]$. Notice that
 \bea
 \xi(r,0)=I_{2m}, ~~~\xi(r,1)=\gamma_{c,h,N}(1) \lb{1.290}
\eea
from $\gamma_{c,h,N}\sim \psi$ with fixed end points for $0 \le r \le r_{1}$. In particular, $\xi(r,\cdot)$ is $C^{1}$ in $t$ for $r\in [r_{1},1]$, hence we can define
\begin{eqnarray}
B_{r}(t)=-J{\dot{\xi}}(r,t)\xi^{-1}(r,t), \hspace*{0.5cm}r_{1} \le r \le 1,~~0 \le t \le 1.\lb{1.110}\end{eqnarray}
Note that $B_{1}(t)={\hat B}_{j}$, $B_{r}(t)$ is continuous in $r,t$.

{\bf Step 2.} Find the homotopy $\zeta(s,h,N)(t_{n})$ connecting $\gamma_{d,h,N}(t_{n})$ to $\xi(r_{1},t_{n})$ with $0 \le s \le r_{1}, 0\le n \le N$.

This task will be divided into five sub-steps.

{\bf Step 2.1.} The discretization of $\xi(r,t)$ with $0 \le r \le r_{1}$

Consider the case the homotopy $\xi(r,\cdot)$ with $0\le r \le r_{1}$. Since $\xi(r,t)$ is continuous in $r,t$, together with the compactness of $[0, r_{1}]\times [0,1]$, then $\xi$ is uniformly continuous and uniformly bounded on $[0, r_{1}]\times [0,1]$, thus $\forall \epsilon >0$, there exists ${\bar h}_{1}>0$ such that for $|r-r'|+|t-t'|\le {\bar h}_{1}$, we have
\bea
|\xi(r,t)-\xi(r',t')|<\epsilon.\lb{1.184}
\eea
And together with the symplecticity of $\xi(r,t)$, there exist two positive constant $c_{3}, c_{4}$ satisfying
\bea
c_{3} \le |\xi(r,t)| \le c_{4}\lb{1.185}
\eea for all $(r,t)\in [0, r_{1}]\times [0,1]$.

While
\bea
|\xi(r,t_{n+1})\xi^{-1}(r,t_{n})-I_{2m}|&=&|[\xi(r,t_{n+1})-\xi(r,t_{n})]\xi^{-1}(r,t_{n})|\nn\\
                                        &\le&|\xi(r,t_{n+1})-\xi(r,t_{n})||\xi^{-1}(r,t_{n})|\nn\\
                                        &\le&\frac{1}{c_{3}}|\xi(r,t_{n+1})-\xi(r,t_{n})|\nn\\
                                        &\le&\frac{1}{c_{3}}\epsilon.\lb{1.186}
\eea
where we use (\ref{1.185}) in the third line above, the last line holds by (\ref{1.184}). Combining (\ref{1.184}) with (\ref{1.186}), there exists ${\bar h}_{1}>0$ sufficiently small such that $\xi(r,t_{n+1})\xi^{-1}(r,t_{n})$ is uniformly close enough to $I_{2m}$ for any $r \in [0,r_{1}]$, where $0 \le n \le N-1$.

{\bf Step 2.2.} The discrete Hamiltonian system and the bilinear form $\mathcal{A}_{r,h,\omega}$, the continuity of ${\tilde B}_{r,n}$ in $r$

Based on the discussion in Step 2.1, by Lemma 2.2, there exists symmetric matrices ${\tilde B}_{r,n}$ which is determined by $\xi(r,t_{n+1})\xi^{-1}(r,t_{n})$.
Since $\xi(r,t_{n+1})\xi^{-1}(r,t_{n})$ is continuous in $r$, it shows that ${\tilde B}_{r,n}$ is continuous in $r$, by Lemma 2.2.
To be precise, suppose that
\begin{eqnarray}
\xi(r,t_{n+1})\xi^{-1}(r,t_{n}):=X_{r,n}=\left(
               \begin{array}{cc}
               X_{r,n,11} & X_{r,n,12}\\
               X_{r,n,21} & X_{r,n,22}\\
               \end{array}
             \right),\lb{1.199}
\end{eqnarray}
where $X_{r,n,11}, X_{r,n,12}, X_{r,n,21}, X_{r,n,22}$ are $m\times m$ continuous matrices in $r$. By Lemma 2.2, we have
\be
{\tilde B }_{r,n}=\left(
              \begin{array}{cc}
                {\tilde A }_{r,n} & {\tilde C }_{r,n} \\
                {\tilde C }_{r,n}^{T} & {\tilde D }_{r,n} \\
              \end{array}
            \right)\lb{1.200}
\ee
with
\bea
{\bar h}_{1}{\tilde C }_{r,n}=(X_{r,n,11}^{-1}-I)^{T},\lb{1.201}
\eea
\bea
{\bar h}_{1}{\tilde A }_{r,n}=X_{r,n,21}(I+{\bar h}_{1}{\tilde C }_{r,n}^{T}),\lb{1.202}
\eea
\bea
{\bar h}_{1}{\tilde D }_{r,n}=-(I+{\bar h}_{1}{\tilde C }_{r,n}^{T})X_{r,n,12}.\lb{1.203}
\eea
It follows that ${\tilde B }_{r,n}$ is continuous in $r$ from (\ref{1.201}), (\ref{1.202}), (\ref{1.203}) and the fact that $X_{r,n}$ is continuous in $r$. Note that $X_{r,n}$ is uniformly bounded from (\ref{1.186}) and (\ref{1.199}), hence ${\tilde B }_{r,n}$ are also uniformly bounded from (\ref{1.201}), (\ref{1.202}) and (\ref{1.203}).

Furthermore we obtain the corresponding discrete Hamiltonian system
\be
\frac{z_{n+1}-z_{n}}{{\bar h}_{1}}=J{\tilde B}_{r,n}\tilde {z}_{n}, \lb{1.140}
\ee
and the bilinear form
\bea
\mathcal{A}_{r,{\bar h}_{1},\omega}(\tilde {z},\tilde {w})=\Sigma_{n=1}^{N}({-J\frac{{z}_{n+1}-z_{n}}{{\bar h}_{1}},\tilde{w}_{n}})-({\tilde B}_{r,n}\tilde{z}_{n}, \tilde{w}_{n})\lb{1.145}
\eea
on ${\tilde W}$.

{\bf Step 2.3.} Well-definedness of $i_{\omega}(\gamma_{r,d,h,N})$

As v) of Remark 2.3 and i) of Remark 2.5, from (\ref{1.140}), denote by
\begin{eqnarray}
\gamma_{r,n+1}=
S_{r,n,{\bar h}_{1}}\gamma_{r,n},\lb{1.222}
\end{eqnarray}
and
\begin{eqnarray}
\gamma_{r,n,s,{\bar h}_{1}}=
S_{r,n,s,{\bar h}_{1}}\gamma_{r,n}.\lb{1.223}
\end{eqnarray}
where $S_{r,n,{\bar h}_{1}}$, $S_{r,n,s,{\bar h}_{1}}$ is given by (\ref{1.213}), (\ref{1.215}) with $A_{r,n}, C_{r,n}, D_{r,n}, F_{r,n}, h$ replacing by ${\tilde A}_{r,n}, {\tilde C}_{r,n}, {\tilde D}_{r,n}, {\tilde C}_{r,n}^{T}, {\bar h}_{1}$ respectively.

By (\ref{1.186}), (\ref{1.201}), (\ref{1.202}) and (\ref{1.203}), there holds
\bea
|{\bar h}_{1}{\tilde B }_{r,n}|\le O(\epsilon), ~~~\forall 0 \le n \le N. \lb{1.207}
\eea
Then it follows that
\bea
|S_{r,n,s,{\bar h}_{1}}-I_{2m}|\le O(\epsilon)\lb{1.208}
\eea
from (\ref{1.207}) and v) of Remark 2.3, i) of Remark 2.5.

Denote by ${\tilde \gamma}_{r,d,{\bar h}_{1},N}$ the fundamental solution to (\ref{1.140}). From (\ref{1.208}), we conclude that ${\tilde \gamma}_{r,c,{\bar h}_{1},N}$ is well defined, hence $i_{\omega}({\tilde\gamma}_{r,c,{\bar h}_{1},N})$ is well defined.

Moreover, if we choose ${h_{1}}=min(h,{\bar h}_{1})$, according to the fact that $\gamma_{c,h,N} \sim \psi$ with fixed end points, note that there hold $i_{\omega}({\tilde\gamma}_{r,d,h_{1},{ N}})=i_{\omega}({\tilde\gamma}_{r,d,{\bar h}_{1},N})=i_{\omega}({\tilde\gamma}_{r,d,h,N})=j$ and $\nu_{\omega}({\tilde\gamma}_{r,d,h_{1}, N}(t_{N}))=\nu_{\omega}({\tilde\gamma}_{r,d,{\bar h}_{1},N}(t_{N}))=\nu_{\omega}({\tilde\gamma}_{r,d,h,N}(t_{N}))=0$ in sense of Corollary 2.8, then the homotopy $\xi(r,t)$ for $h_{1}$ is the same as one for $h$, therefore the discretization process above remains true, finally we can take $h_{1}$ as $h$.

{\bf Step 2.4.} Find $\zeta(s,h,N)(t_{n})$ with $i_{\omega}(\zeta(s,h,N))=j$, $\nu_{\omega}(\zeta(s,h,N)(t_{N}))=0$ for $0 \le s \le r_{1}$.

Using Lemma 2.2, we have
\be
{\tilde\gamma}_{r,d,h,N}(t_{n+1})=\xi(r,t_{n+1})\xi^{-1}(r,t_{n}){\tilde\gamma}_{r,d,h,N}(t_{n}),\lb{1.141}
\ee
Note that since ${\tilde\gamma}_{r,d,h,N}(0)=\xi(r,0)=I_{2m}$, there holds
\bea
{\tilde\gamma}_{r,d,h,N}(t_{n})=\xi(r,t_{n}).\lb{1.142}
\eea
We can define $\zeta(s,h,N)(t_{n})$ by
\bea
\zeta(s,h,N)(t_{n})={\tilde\gamma}_{r,d,h,N}(t_{n}),\lb{1.187}
\eea
where $s=r$. Since $\xi(r,t)$ is continuous in $r$, then so is $\xi(r,t_{n})$, thus so is ${\tilde\gamma}_{r,d,h,N}(t_{n})$ from (\ref{1.142}), it is followed that  $\zeta(s,h,N)(t_{n})$ is continuous in $s$.

By (\ref{1.142}) and (\ref{1.187}), Definition 2.6, Corollary 2.8, Lemma 2.10 and Lemma 4.4, together with $\xi(r,t_{N})=\gamma_{c,h,N}(t_{N})$ from Step 1 and $\nu_{\omega}(\gamma_{d,h,N}(t_{N}))=0$, we conclude that $i_{\omega}(\zeta(s,h,N))=j$, $\nu_{\omega}(\zeta(s,h,N)(t_{N}))=0$ for $0 \le s \le r_{1}$.

{\bf Step 2.5.} $\zeta(s,h,N)(t_{n})$ with $\zeta(0,h,N)(t_{n})=\gamma_{d,h,N}(t_{n})$ and $\zeta(r_{1},h,N)(t_{n})=\xi(r_{1},t_{n})$

By (\ref{1.288}), (\ref{1.142}) and (\ref{1.187}), the results hold.\\

{\bf Step 3.} Find the homotopy $\zeta(s,h,N)(t_{n})$ connecting $\xi(r_{1},t_{n})$ to $\gamma_{r_{1},d,h,N}(t_{n})$ with $r_{1} \le s \le 1+r_{1} , 0\le n \le N$.

In general, there is a gap between ${\tilde\gamma}_{r_{1},d,h,N}(t_{n})$ and $\gamma_{r_{1},d,h,N}(t_{n})$, that is ${\tilde\gamma}_{r_{1},d,h,N}(t_{n})\ne \gamma_{r_{1},d,h,N}(t_{n})$ for $1 \le n \le N$, so the goal of this step is to find the homotopy $\zeta(s,h,N)(t_{n})$ which connects ${\tilde\gamma}_{r_{1},d,h,N}(t_{n})$ to $\gamma_{r_{1},d,h,N}(t_{n})$, where ${\gamma}_{r_{1},d,h,N}$ is the fundamental solution to (\ref{1.234}) with ${\bar h}_{2}$ replacing by $h$ and $r=r_{1}$. It will be realized in seven sub-steps.

{\bf Step 3.1.} The gap between ${\tilde\gamma}_{r_{1},d,h,N}(t_{n})$ and $\gamma_{r_{1},d,h,N}(t_{n})$

Since $\xi(r,t)$ is uniformly bounded on $[0,r_{1}]\times[0,1]$, ${\tilde\gamma}_{r,d,h,N}(t_{n})=\xi(r,t_{n})$,  together with (\ref{1.185}) in Step 2.1, it yields that
\bea
c_{3} \le |{\tilde\gamma}_{r_{1},d,h,N}(t_{n})| \le c_{4}\lb{1.261}
\eea
for $0 \le n \le N$.

According to A.1, $\forall \epsilon>0$, there exists $h>0$ sufficiently small such that
\bea
|{\tilde\gamma}_{r_{1},d,h,N}(t_{n})-\gamma_{r_{1},d,h,N}(t_{n})|\le \epsilon.\lb{1.262}
\eea
Therefore, we have
\bea
&&|\gamma_{r_{1},d,h,N}(t_{n}){\tilde\gamma}_{r_{1},d,h,N}(t_{n})^{-1}-I_{2m}|\nn\\
&=&|\gamma_{r_{1},d,h,N}(t_{n})-{\tilde\gamma}_{r_{1},d,h,N}(t_{n})||{\tilde\gamma}_{r_{1},d,h,N}(t_{n})^{-1}|\nn\\
&\le&\frac{\epsilon}{c_{3}},\lb{1.263}
\eea
where we use (\ref{1.261}) and (\ref{1.262}) in last line.

{\bf Step 3.2.} The continuity of ${\bar B}_{r,n}$ in $r$, the discrete Hamiltonian system and the bilinear form $\mathcal{A}_{r,h,\omega}$

Set
\begin{eqnarray}
\gamma_{r_{1},d,h,N}(t_{n}){\tilde\gamma}_{r_{1},d,h,N}(t_{n})^{-1}:={\bar X}_{n}=\left(
               \begin{array}{cc}
               {\bar X}_{n,11} & {\bar X}_{n,12}\\
               {\bar X}_{n,21} & {\bar X}_{n,22}\\
               \end{array}
             \right),\lb{1.264}
\end{eqnarray}
where ${\bar X}_{n,11}, {\bar X}_{n,12}, {\bar X}_{n,21}, {\bar X}_{n,22}$ are $m\times m$  matrices. By Lemma 2.2 and iv) of Remark 2.3, we have
\be
{\bar B }_{n}=\left(
              \begin{array}{cc}
                {\bar A }_{n} & {\bar C }_{n} \\
                {\bar C }_{n}^{T} & {\bar D }_{n} \\
              \end{array}
            \right)\lb{1.265}
\ee
with
\bea
h{\bar C }_{n}=({\bar X}_{n,11}^{-1}-I)^{T},\lb{1.266}
\eea
\bea
h{\bar A }_{n}={\bar X}_{n,21}(I+h{\bar C }_{n}^{T}),\lb{1.267}
\eea
\bea
h{\bar D }_{n}=-(I+h{\bar C }_{n}^{T}){\bar X}_{n,12}.\lb{1.268}
\eea
Moreover, by (\ref{1.263}), we have
\bea
h{\bar C }_{n}=O(\epsilon),\lb{1.277}
\eea
\bea
h{\bar A }_{n}=O(\epsilon),\lb{1.278}
\eea
\bea
h{\bar D }_{n}=O(\epsilon).\lb{1.279}
\eea

By the similar idea to (\ref{1.33}) in section 2.2, let
\bea
{\bar B}_{r,n}=\left(
              \begin{array}{cc}
                {\bar A}_{r,n} & {\bar C}_{r,n} \\
                {\bar C}_{r,n}^{T} & {\bar D}_{r,n} \\
              \end{array}
            \right):=r{\bar B}_{n}.\lb{1.269}
\eea
where $0 \le r \le 1$. It is obvious that ${\bar B}_{r,n}$ is continuous in $r$.
Moreover from (\ref{1.269}), it yields the discrete Hamiltonian system
\be
\frac{z_{n+1}-z_{n}}{h}=J{\bar B}_{r,n}\tilde {z}_{n},\lb{1.238}
\ee
and the bilinear form
\bea
\mathcal{A}_{r,h,\omega}(\tilde {z},\tilde {w})=\Sigma_{n=1}^{N_{2}}({-J\frac{{z}_{n+1}-z_{n}}{h},\tilde{w}_{n}})-({\bar B}_{r,n}\tilde{z}_{n} \tilde{w}_{n}),\lb{1.239}
\eea
on ${\tilde W}$.

{\bf Step 3.3.} The continuous connecting orbits ${\bar\gamma}_{r,d,h,n}$ between ${\tilde\gamma}_{r_{1},d,h,N}(t_{n})$ and $\gamma_{r_{1},d,h,N}(t_{n})$

For consistency of notations, set
\bea
{\bar S}_{n,h}={\bar X}_{n}=\left(
\begin{array}{cc}
(I+h{\bar C}_{n}^{T})^{-1} & -h(I+h{\bar C}_{n}^{T})^{-1}{\bar D}_{n}\\
h{\bar A}_{n}(I+h{\bar C}_{n}^{T})^{-1} & -h^{2}{\bar A}_{n}(I+h{\bar C}_{n}^{T})^{-1}{\bar D}_{n}+(I+h{\bar C}_{n})\\
\end{array}
\right).\lb{1.273}
\eea
Note that from (\ref{1.264}), we have
\bea
\gamma_{r_{1},d,h,N}(t_{n})={\bar S}_{n,h}{\tilde\gamma}_{r_{1},d,h,N}(t_{n}).\lb{1.280}
\eea
By (\ref{1.263}), (\ref{1.264}), there holds
\bea
|{\bar S}_{n,h}-I_{2m}|\le O(\epsilon).\lb{1.275}
\eea
As in section 2, we obtain
\bea
\eta_{r,n,h}={\bar S}_{r,n,h}\eta_{0,n,h}\lb{1.270}
\eea
from (\ref{1.269}) and (\ref{1.238}), where $\eta_{0,n,h},\eta_{r,n,h}$ possess the same meaning as $\gamma_{n}, \gamma_{n,s,h}$ in (\ref{1.39}),
\begin{eqnarray}
{\bar S}_{r,n,h}=
\left(
\begin{array}{cc}
(I+rh{\bar C}_{n}^{T})^{-1} & -rh(I+rh{\bar C}_{n}^{T})^{-1}{\bar D}_{n}\\
rh{\bar A}_{n}(I+rh{\bar C}_{n}^{T})^{-1} & -r^{2}h^{2}{\bar A}_{n}(I+rh{\bar C}_{n}^{T})^{-1}{\bar D}_{n}+(I+rh{\bar C}_{n})\\
\end{array}
\right).\lb{1.271}
\end{eqnarray}
and (\ref{1.270}) is the counterpart of (\ref{1.39}) in Remark 2.5.

It is followed that ${\bar S}_{1,n,h}=\gamma_{r_{1},d,h,N}(t_{n}){\tilde\gamma}_{r_{1},d,h,N}(t_{n})^{-1}$ from Lemma 2.2 and Lemma 2.4. Moreover, we have
\bea
{\bar S}_{r,n,h}=\left\{\begin{array}{ll}
      I_{2m}, \hspace*{0.2cm}r=0. \\
      {\bar S}_{n,h},\hspace*{0.3cm} r=1. \end{array}\right.\lb{1.272}
\eea
Denote by ${\bar\gamma}_{r,d,h,n}$ the solution to the following equation with coefficient matrix ${\bar S}_{r,n,h}$ and the initial value ${\bar\gamma}_{0,d,h,n}={\tilde\gamma}_{r_{1},d,h,N}(t_{n})$, that is
\bea
{\bar\gamma}_{r,d,h,n}={\bar S}_{r,n,h}{\tilde\gamma}_{r_{1},d,h,N}(t_{n}).\lb{1.274}
\eea
From (\ref{1.280}) and (\ref{1.272}), it forces that
${\bar\gamma}_{1,d,h,n}=\gamma_{r_{1},d,h,N}(t_{n})$. Thus ${\bar\gamma}_{r,d,h,n}$ is the connecting orbit between ${\tilde\gamma}_{r_{1},d,h,N}(t_{n})$ and $\gamma_{r_{1},d,h,N}(t_{n})$.

Furthermore, when $n=0$, then there holds
\bea
{\tilde\gamma}_{r_{1},d,h,N}(t_{0})=\gamma_{r_{1},d,h,N}(t_{0})=I_{2m}, \lb{1.281}
\eea
hence from ({\ref{1.264}}) and ({\ref{1.273}}), it follows that
\bea
{\bar S}_{0,h}={\bar X}_{0}=I_{2m},\lb{1.282}
\eea
By ({\ref{1.266}}), ({\ref{1.267}}) and ({\ref{1.268}}), it shows that
\bea
{\bar A}_{0}={\bar C}_{0}={\bar D}_{0}=O_{2m}.\lb{1.283}
\eea
Then substituting ${\bar A}_{0}={\bar C}_{0}={\bar D}_{0}$ into (\ref{1.271}), we have
\bea
{\bar S}_{r,n,h}=I_{2m}.\lb{1.284}
\eea
According to (\ref{1.274}), it leads to
\bea
{\bar\gamma}_{r,d,h,0}=I_{2m}.\lb{1.285}
\eea
for $\forall r \in [0,1].$

{\bf Step 3.4.} The continuity and non-degenerate of ${\bar\gamma}_{r,d,h,N}$

Since ${\bar S}_{r,n,h}$ is continuous in $r$, so is ${\bar\gamma}_{r,d,h,n}$. On the other hand, we have
\begin{eqnarray}
&&|{\bar\gamma}_{r,d,h,n}-{\tilde\gamma}_{r_{1},d,h,N}(t_{n})|\nn\\
&=&|{\bar S}_{r,n,h}{\tilde\gamma}_{r_{1},d,h,N}(t_{n})-{\tilde\gamma}_{r_{1},d,h,N}(t_{n})|\nn\\
&=&|{\bar S}_{r,n,h}-I_{2m}||{\tilde\gamma}_{r_{1},d,h,N}(t_{n})|\nn\\
&\le&O(\epsilon).\lb{1.296}
\end{eqnarray}
The third equality holds from (\ref{1.261}), (\ref{1.266}), (\ref{1.267}), (\ref{1.268}), (\ref{1.275}) and (\ref{1.271}).

Since $\nu_{\omega}({\tilde\gamma}_{r_{1},d,h,N}(t_{N}))=0$, then there exists a neighborhood of ${\tilde\gamma}_{r_{1},d,h,N}(t_{N})$ in $Sp^{*}_{\omega}(2m)$, hence there is a positive distance between ${\tilde\gamma}_{r_{1},d,h,N}(t_{N})$ and $Sp^{0}_{\omega}(2m)$. Consequently, we conclude that $\nu_{\omega}({\bar\gamma}_{r,d,h,N})=0$ from (\ref{1.296}).

{\bf Step 3.5.} Well-posedness of $i_{\omega}({\bar\gamma}_{r,d,h,n})$

Since
\begin{eqnarray}
&&|{\bar\gamma}_{r,d,h,n+1}-{\bar\gamma}_{r,d,h,n}|\nn\\
&=&|{\bar S}_{r,n+1,h}{\tilde\gamma}_{r_{1},d,h,N}(t_{n+1})-{\bar S}_{r,n,h}{\tilde\gamma}_{r_{1},d,h,N}(t_{n})|\nn\\
&=&|{\bar S}_{r,n+1,h}{\tilde\gamma}_{r_{1},d,h,N}(t_{n+1})-{\bar S}_{r,n+1,h}{\tilde\gamma}_{r_{1},d,h,N}(t_{n})+{\bar S}_{r,n+1,h}{\tilde\gamma}_{r_{1},d,h,N}(t_{n})-{\bar S}_{r,n,h}{\tilde\gamma}_{r_{1},d,h,N}(t_{n})|\nn\\
&\le&|{\bar S}_{r,n+1,h}{\tilde\gamma}_{r_{1},d,h,N}(t_{n+1})-{\bar S}_{r,n+1,h}{\tilde\gamma}_{r_{1},d,h,N}(t_{n})|+|{\bar S}_{r,n+1,h}{\tilde\gamma}_{r_{1},d,h,N}(t_{n})-{\bar S}_{r,n,h}{\tilde\gamma}_{r_{1},d,h,N}(t_{n})|\nn\\
&=&|{\bar S}_{r,n+1,h}||{\tilde\gamma}_{r_{1},d,h,N}(t_{n+1})-{\tilde\gamma}_{r_{1},d,h,N}(t_{n})|+|{\bar S}_{r,n+1,h}-{\bar S}_{r,n,h}||{\tilde\gamma}_{r_{1},d,h,N}(t_{n})|\nn\\
&=&|{\bar S}_{r,n+1,h}||{\tilde\gamma}_{r_{1},d,h,N}(t_{n+1}){\tilde\gamma}_{r_{1},d,h,N}(t_{n})^{-1}-I_{2m}||{\tilde\gamma}_{r_{1},d,h,N}(t_{n})|+|{\bar S}_{r,n+1,h}-{\bar S}_{r,n,h}||{\tilde\gamma}_{r_{1},d,h,N}(t_{n})|\nn\\
&\le&O(\epsilon).\lb{1.276}
\end{eqnarray}
where the first equality uses (\ref{1.274}), the last inequality is followed from (\ref{1.261}), (\ref{1.263}), (\ref{1.277}), (\ref{1.278}), (\ref{1.279}) and (\ref{1.296}).

Thus $i_{\omega}({\bar\gamma}_{r,d,h,n})$ is well defined by Lemma 2.2.

{\bf Step 3.6.} Find $\zeta(s,h,N)(t_{n})$ satisfying $i_{\omega}(\zeta(s,h,N))=j$ and $\nu_{\omega}(\zeta(s,h,N)(t_{N}))=0$ with $r_{1}\le s \le 1+r_{1}$.

We can define $\zeta(s,h,N)(t_{n})$ by
\bea
\zeta(s,h,N)(t_{n})={\bar\gamma}_{r,d,h,n},\lb{1.252}
\eea
where $s=r+r_{1}, 0 \le r \le 1$.
Since ${\bar\gamma}_{r,d,h,n}$ is continuous in $r$, so $\zeta(s,h,N)$ is continuous in $s$. Moreover, by (\ref{1.285}), (\ref{1.296}) , (\ref{1.276}) and (\ref{1.252}), Definition 2.6 , Corollary 2.8, Lemma 2.10 and Lemma 4.4, together with $i_{\omega}({\bar\gamma}_{0,d,h,N})=i_{\omega}(\psi)=j$ and $\nu_{\omega}({\bar\gamma}_{0,d,h,N})=0$, we conclude that $i_{\omega}(\zeta(s,h,N))=j$, $\nu_{\omega}(\zeta(s,h,N)(t_{N}))=0$ for $r_{1}\le s \le 1+r_{1}$.

{\bf Step 3.7.} $\zeta(r_{1},h,N)(t_{n})=\xi(r_{1},t_{n}), \zeta(1+r_{1},h,N)(t_{n})=\gamma_{r_{1},d,h,N}(t_{n})$

By (\ref{1.280}), (\ref{1.272}) and (\ref{1.274}), there hold ${\bar\gamma}_{0,d,h,n}={\tilde\gamma}_{r_{1},d,h,N}(t_{n}),{\bar\gamma}_{1,d,h,n}=\gamma_{r_{1},d,h,N}(t_{n})$, together with (\ref{1.142}), it follows that ${\bar\gamma}_{0,d,h,n}=\xi(r_{1},t_{n})$. In addition to (\ref{1.252}), the results are followed.\\

{\bf Step 4.} Find the homotopy $\zeta(s,h,N)(t_{n})$ connecting ${\gamma}_{r_{1},h,d,N}$ to ${\hat{\gamma}}_{j,h,d,N}$ with $1+r_{1} \le s \le 2 , 0\le n \le N$.

This scheme  consists of four sub-steps.

{\bf Step 4.1.} The discrete Hamiltonian system, the bilinear form $\mathcal{A}_{r,h,\omega}$ and the continuity of ${\tilde B}_{r,n}$ in $r$

Since the homotopy $\xi(r,\cdot)$ is $C^{1}$ in $t$ for $r_{1}\le r \le 1$, we can obtain symmetric matrices $B_{r}(t)$ given by (\ref{1.110}) and the corresponding continuous Hamiltonian system
\be
{\dot z}=JB_{r}(t)z.\lb{1.143}
\ee
Note that $B_{r}(t)$ is continuous in $r,t$. From (\ref{1.143}), it yields the discrete Hamiltonian system
\be
\frac{z_{n+1}-z_{n}}{h}=JB_{r,n}\tilde {z}_{n},\lb{1.144}
\ee
and the bilinear form
\bea
\mathcal{A}_{r,h,\omega}(\tilde {z},\tilde {w})=\Sigma_{n=1}^{N_{2}}({-J\frac{{z}_{n+1}-z_{n}}{h},\tilde{w}_{n}})-(B_{r,n}\tilde{z}_{n} \tilde{w}_{n}),\lb{1.146}
\eea
on ${\tilde W}$, where $B_{r,n}=B_{r}(t_{n})$. Since $B_{r}(t)$ is continuous in $r$, so is $B_{r,n}$. Suppose $B_{r,n}$ has the form
with
\be
B_{r,n}=\left(
              \begin{array}{cc}
                A_{r,n} & C_{r,n} \\
                C_{r,n}^{T} & D_{r,n} \\
              \end{array}
            \right).\lb{1.224}
\ee

{\bf Step 4.2.} Well-posedness of $i_{\omega}(\gamma_{r,d,h,N})$ and the continuity of $\gamma_{r,d, h_{2},N}$ in $r$

Since $B_{r}(t)$ is continuous in $r,t$, together with the compactness of $[r_{1},1]\times [0,1]$, then $B_{r}(t)$ is uniformly bounded on $[r_{1},1]\times [0,1]$, that is there exists positive constant $c_{5}$ satisfying $|B_{r}(t)|\le c_{5}$ for all $(r,t) \in [r_{1},1]\times [0,1]$. Therefore, $\forall \epsilon>0$, there exists ${\bar h}_{2}>0$ sufficiently small such that
\bea
|{\bar h}_{2}B_{r}(t)|\le \epsilon \lb{1.225}
\eea
for all $(r,t) \in [r_{1},1]\times [0,1]$. In particular, we have
\bea
|{\bar h}_{2}B_{r,n}|\le \epsilon \lb{1.226}
\eea
for $0 \le n \le N$. The discrete Hamiltonian system with the coefficient matrix $B_{r,n}$ is
\be
\frac{z_{n+1}-z_{n}}{{\bar h}_{2}}=JB_{r,n}\tilde {z}_{n}.\lb{1.234}
\ee

By Lemma 2.2 and v) of Remark 2.3 and i) of Remark 2.5, it yields that
\begin{eqnarray}
S_{r,n,{\bar h}_{2}}=
\left(
              \begin{array}{cc}
                (I+{\bar h}_{2}C_{r,n}^{T})^{-1} & -h(I+{\bar h}_{2}C_{r,n}^{T})^{-1}D_{r,n} \\
                {\bar h}_{2}A_{r,n}(I+{\bar h}_{2}C_{r,n}^{T})^{-1} & -{\bar h}_{2}^{2}A_{r,n}(I+{\bar h}_{2}C_{r,n}^{T})^{-1}D_{r,n}+I+{\bar h}_{2}C_{r,n} \\
              \end{array}
            \right)\lb{1.227}
            \end{eqnarray}
satisfying
\be
\gamma_{r,n+1}=S_{r,n,{\bar h}_{2}}\gamma_{r,n}.\lb{1.228}
\ee
and
\begin{eqnarray}
S_{r,n,s',{\bar h}_{2}}=
\left(
\begin{array}{cc}
(I+s'{\bar h}_{2}C_{r,n}^{T})^{-1} & -s'h(I+s'{\bar h}_{2}C_{r,n}^{T})^{-1}D_{n}\\
s'{\bar h}_{2}A_{r,n}(I+s'{\bar h}_{2}C_{r,n}^{T})^{-1} & -{s'}^{2}{\bar h}_{2}^{2}A_{r,n}(I+s'{\bar h}_{2}C_{r,n}^{T})^{-1}D_{r,n}+(I+s'{\bar h}_{2}C_{r,n})\\
\end{array}
\right),\lb{1.229}
\end{eqnarray}
\bea
S_{r,n,s',{\bar h}_{2}}=\left\{\begin{array}{ll}
 I_{2m}, \hspace*{0.2cm}s'=0. \cr
 S_{r,n,{\bar h}_{2}},\hspace*{0.3cm} s'=1. \end{array}\right.\lb{1.230}
\eea
satisfying
\begin{eqnarray}
\gamma_{r,n,s',{\bar h}_{2}}=
S_{r,n,s',{\bar h}_{2}}\gamma_{r,n},\lb{1.231}
\end{eqnarray}
with $\gamma_{r,n,0,{\bar h}_{2}}=\gamma_{r,n}, ~~\gamma_{r,n,1,{\bar h}_{2}}=\gamma_{r,n+1}=\gamma_{r,n+1,0,{\bar h}_{2}}$ and $ \gamma_{r,n+1,1,{\bar h}_{2}}=\gamma_{r,n+2}$, where $0 \le s' \le 1$. Note that $S_{r,n,{\bar h}_{2}}$ is continuous in $r$ from the fact that $B_{r,n}$ is continuous in $r$.

By (\ref{1.226}) and (\ref{1.227}), there holds
\begin{eqnarray}
|S_{r,n,{\bar h}_{2}}-I_{2m}| \le O(\epsilon),\lb{1.232}
\end{eqnarray}
Hence from (\ref{1.229}), we have
\begin{eqnarray}
|S_{r,n,s',{\bar h}_{2}}-I_{2m}| \le O(\epsilon).\lb{1.233}
\end{eqnarray}

Denote by $\gamma_{r,d,{\bar h}_{2},N}$ the fundamental solution to (\ref{1.234}), in addition to (\ref{1.228}), there holds
\be
\gamma_{r,d,{\bar h}_{2},N}(t_{n+1})=S_{r,n,{\bar h}_{2}}\gamma_{r,d,{\bar h}_{2},N}(t_{n}).\lb{1.235}
\ee

Since $\gamma_{r,d,{\bar h}_{2},N}(t_{0})=I_{2m}$, it follows that $\gamma_{r,d,{\bar h}_{2},N}$ is continuous in $r$ from iii) of A.1. Moreover from (\ref{1.233}), we conclude that $\gamma_{r,c,{\bar h}_{2},N}$ is well defined, hence $i_{\omega}(\gamma_{r,c,{\bar h}_{2},N})$ is well defined.

Moreover, if we choose ${h_{2}}=min(h,{\bar h}_{2})$, according to the fact that $\psi \sim {\hat\beta}_{j}$, note that there hold $i_{\omega}(\gamma_{r,d,h_{2},{ N}})=i_{\omega}(\gamma_{r,d,{\bar h}_{2},N})=i_{\omega}(\gamma_{r,d,h,N})=j$ and $\nu_{\omega}(\gamma_{r,d,h_{2}, N}(t_{N}))=\nu_{\omega}(\gamma_{r,d,{\bar h}_{2},N}(t_{N}))=\nu_{\omega}(\gamma_{r,d,h,N}(t_{N}))=0$ in sense of Corollary 2.8, then the homotopy $\xi(r,t)$ for $h_{2}$ is the same as one for $h$, therefore the discretization process above remains true, finally we can take $h_{2}$ as $h$.

{\bf Step 4.3.} Find $\zeta(s,h,N)(t_{n})$ satisfying $i_{\omega}(\zeta(s,h,N))=j$, $\nu_{\omega}(\zeta(s,h,N)(t_{N}))=0$ for $1+ r_{1} \le s \le 2$.

We can define $\zeta(s,h,N)(t_{n})$ by
\bea
\zeta(s,h,N)(t_{n})=\gamma_{r,d,h,N}(t_{n}),\lb{1.188}
\eea
where $s=1+r, r_{1} \le r \le 1$.

Since $\gamma_{r,d,h,N}$ is continuous in $r$, $\zeta(s,h,N)$ is continuous with respect to $s$. By $\xi(r,\cdot) \in \hat{\mathcal P}_{1,\omega}^{*}(2m)$, then there exists a neighborhood of $\xi(r,t_{N})$ in $Sp^{*}_{\omega}(2m)$, hence there is a positive distance $\delta(r)$ between $\xi(r,t_{N})$ and $Sp^{0}_{\omega}(2m)$. Since $\xi(r,t_{N})$ is continuous in $r$, we can take
 \bea
 \delta_{0}=min_{r_{1} \le r \le 1}\delta(r)>0.\lb{1.183}
 \eea
While, by ii) of A.1, for and $\epsilon >0, $ there exist ${\bar h}_{3}>0$ sufficiently small such that for $|t-t_{n}|\le{\bar h}_{3},$ we have
\bea
|\xi(r,t)-\gamma_{r,d,{\bar h}_{3},N}(t_{n})|<\epsilon.\lb{1.236}
\eea
Furthermore it follows that $\nu_{\omega}(\gamma_{r,d,{\bar h}_{3},N}(t_{N}))=0$ and $\gamma_{r,d,{\bar h}_{3},N}(t_{N})$ is in the same connectedness component as that of $\xi(r,t_{N})$ from A.1 and (\ref{1.183}), (\ref{1.236}). In addition to the fact that $\xi(r_{1},t)=\psi$ and $\xi(1,t)={\hat\beta}_{j}$, it is forced that $i_{\omega}(\zeta(s,{\bar h}_{3},N))=i_{\omega}(\gamma_{r,d,{\bar h}_{3},N})=j$, $\nu_{\omega}(\zeta(s,{\bar h}_{3},N)(t_{N}))=\nu_{\omega}(\gamma_{r,d,{\bar h}_{3},N}(t_{N}))=0$ for $r_{1} \le r \le 1$ from Definition 2.6, , Corollary 2.8, Lemma 2.10, Lemma 4.4 and (\ref{1.188}).

Summarizing, based on the discussion above and the fact that $\psi \sim {\hat \beta}_{j}$, we can choose $h_{2}=min(h,{\bar h}_{2},{\bar h}_{3})$ and obtain the same results. Note that since there hold $i_{\omega}(\gamma_{r,d,h,{\bar N}})=i_{\omega}(\gamma_{r,d,h,N})=j$ and $\nu_{\omega}(\gamma_{r,d,h,{\bar N}}(1))=\nu_{\omega}(\gamma_{r,d,h,N}(1))=0$ in sense of Corollary 2.7, then the homotopy $\xi(r,t)$ for $\bar N$ is the same as one for $N$, therefore the discretization process above remains true, specifically $\delta$ in (\ref{1.183}) keeps unchanged, finally we can take $\bar N$ as $N$. Note that according to A.1, it may be concluded that the larger $N$ is, the less $|\xi(r,1)-\gamma_{r,d,h,N}(1)|$ is. Thus for $0 \le n \le N$,  we conclude that end matrices $\xi(r,1)$ of a family of the fundamental solutions of  the corresponding discrete Hamiltonian system (\ref{1.4}) with the symmetric matrices $B_{r,n}$ which corresponds to the homotopy $\xi(r,\cdot)$ with $0\le r < r_{1}$  and end matrices of a family of the fundamental solutions of  the corresponding discrete Hamiltonian system (\ref{1.4}) with the symmetric matrices $B_{r,n}$ which corresponds to $B_{r}(t)$ given by (\ref{1.110}) for $r_{1} \le r \le 1$ are in the same connectedness component as that of ${\hat\beta}_{j}(1)$, therefore their Maslov-type index are also $j$ and end matrices are in $Sp(2n)_{\omega}^*$, that is $i_{\omega}(\zeta(s,h,N))=j$, $\nu_{\omega}(\zeta(s,h,N)(t_{N}))=0$ for $1+ r_{1} \le s \le 2$.

{\bf Step 4.4.} $\zeta(1+r_{1},h,N)(t_{n})={\gamma}_{r_{1},h,d,n}, \zeta(2,h,N)(t_{n})={\hat {\gamma}}_{j,h,d,N}(t_{n})$

Note that ${\gamma}_{1,h,d,N}$ is the fundamental solution to (\ref{1.234}) with $r=1$ and ${\bar h}_{2}$ replacing by $h$, together with (\ref{1.188}), we have $\zeta(2,h,N)(t_{n})={\gamma}_{1,h,d,N}(t_{n})$. Since ${\hat{\gamma}}_{j,h,d,N}$ is the fundamental solution to the discrete Hamiltonian system with coefficient matrix ${\hat B}_{j,n}$, in addition to $B_{1,n}={\hat B}_{j,n}$, it forces that $\zeta(2,h,N)(t_{n})={\gamma}_{1,h,d,N}(t_{n})={\hat {\gamma}}_{j,h,d,N}(t_{n})$. By (\ref{1.188}), there exists $\zeta(1+r_{1},h,N)(t_{n})= {\gamma}_{r_{1},h,d,n}$.\hfill\hb

In sum, from the discussions above, we can take the same $h$ in the Step 2, 3 and 4, thus obtain the following homotopy $\zeta(s,h,N)$ connecting ${\gamma}_{d,h,N}$ to ${\hat{\gamma}}_{r,h,d,N}(t_{n})$ for the discrete case
\be \zeta(s,h,N)(t_{n}) = \left\{\begin{array}{lll}
                    {\tilde {\gamma}}_{r,h,d,N}(t_{n}), & s=r ~~\text{and} ~~0 \le r \le r_{1}. \\
                    {\bar {\gamma}}_{r,h,d,n}, & s=r+r_{1} ~~\text{and} ~~0 \le r \le 1.  \\
                    {\gamma}_{r,h,d,N}(t_{n}), & s=r+1 ~~\text{and} ~~r_{1} \le r \le 1.  \\
                      \end{array}
\right.\lb{1.287}  \ee
where $0 \le s \le 2$, $0 \le n \le N$, ${\tilde {\gamma}}_{r,h,d,N}$ is the fundamental solution to (\ref{1.140}) with $h_{1}$ replacing by $h$, ${\bar {\gamma}}_{r,h,d,n}$ is given by (\ref{1.274}), and ${\gamma}_{r,h,d,N}$ is the fundamental solution to (\ref{1.234}) with ${\bar h}_{2}$ replacing by $h$.

\setcounter{equation}{0}
\section{Proof of the main results}

{\bf Proof of Theorem 1.1} Suppose $i_{\omega}(\gamma_{c,h,N})=j$, the proof  is divided into two steps.

{\bf Step 1.} $\nu_{\omega}(\gamma_{d,h,N}(1))=0.$

 Based on the discussions in Section 5 and Lemma 4.4, the proof of this theorem can be reduced to the study of standard paths ${\hat {\gamma}}_{j,h,d,N}$ and the corresponding matrices ${\hat B}_{j,n}$, where ${\hat B}_{j,n}={\hat B}_{j}(t_{n})$,  ${\hat B}_{j}$ are defined by (\ref{1.104}),(\ref{1.106}),(\ref{1.108}). The proof will fall into three cases.

{\bf Case 1.} $\omega =1$.\\
{\bf Case 1.1.} $n=1$ and $j\notin 2{\bf Z}\setminus \{0\}$ or $n\geq 2$.
The corresponding ${\hat B}_{j}$ to standard paths ${\hat\beta}_{j}$ is the constant matrix, thus ${\hat B}_{j,n}\equiv {\hat B}_{j}$ for $n=1,\dots, N$. We only prove the case $n \geq 2$ and $(-1)^{m+j}=1$, the rest cases can be dealt with similarly.

In this case, the bilinear form (\ref{1.47}) is
\bea
\mathcal{A}_{h,1}(\tilde {z},\tilde {w})=\Sigma_{n=1}^{N}({-J\frac{{z}_{n+1}-z_{n}}{h},\tilde{w}_{n}})-({\hat B}_{j}\tilde{z}_{n}. \tilde{w}_{n}),\lb{1.55}
\eea
By (\ref{1.27}),(\ref{1.29}),(\ref{1.30}),(\ref{1.31}),(\ref{1.32}),(\ref{1.36}),(\ref{1.47}), the action space ${\tilde W}$ of $\mathcal{A}_{h,1}$ can be written as
\bea
{\tilde W}=E_{0}\oplus\bigg(\oplus_{k=1}^{[\frac{N-1}{2}]}\bigg[ E_{k}\oplus E_{-k}\bigg]\bigg), ~~~\text{if $N$ is odd.}\lb{1.70}\eea
or
\bea
{\tilde W}=E_{0}\oplus\bigg(\oplus_{k=1}^{[\frac{N}{2}]-1}\bigg[ E_{k}\oplus E_{-k}\bigg]\bigg)\oplus \bigg(E_{\frac{N}{2}}\oplus E_{-\frac{N}{2}}\bigg),~~~\text{if $N$ is even.}\lb{1.71}\eea
Where the invariant subspaces
\begin{eqnarray}
E_{0}={\mathcal{G}}, \mathcal{G}\in {\bf R}^{2m},\nn
\end{eqnarray}
\begin{eqnarray}
E_{k}=\oplus_{n=1}^{N}\bigg[\left(
\begin{array}{cc}
sin(2n+1)\theta_{k}\hspace*{1.5mm}a_{k}\\
cos2n\theta_{k}\hspace*{1.5mm}a_{k}
\end{array}
\right)
\oplus
\left(
\begin{array}{cc}
-cos(2n+1)\theta_{k}\hspace*{1.5mm}b_{k}\\
sin2n\theta_{k}\hspace*{1.5mm}b_{k}
\end{array}
\right)\bigg],\nn
\end{eqnarray}
\begin{eqnarray}
E_{-k}=\oplus_{n=1}^{N}\bigg[\left(
\begin{array}{cc}
-sin(2n+1)\theta_{k}\hspace*{1.5mm}c_{k}\\
cos2n\theta_{k}\hspace*{1.5mm}c_{k}
\end{array}
\right)
\oplus
\left(
\begin{array}{cc}
cos(2n+1)\theta_{k}\hspace*{1.5mm}d_{k}\\
sin2n\theta_{k}\hspace*{1.5mm}d_{k}
\end{array}
\right)\bigg],\nn
\end{eqnarray}
are the eigensubspaces associated to the eigenvalues $\lambda_{0}=0$, $\lambda_{k}=\frac{2sin\theta_{k}}{h}$ and $\lambda_{-k}=-\frac{2sin\theta_{k}}{h}$ respectively, $a_{k}, b_{k}, c_{k}, d_{k} \in {\bf R}^{m}, h=\frac{1}{N}$ and $\theta_{k}=\frac{k\pi}{N}$.

Also the invariant subspaces
\begin{eqnarray}
E_{\frac{N}{2}}=\oplus_{n=1}^{N}\bigg[(-1)^{n}\left(
\begin{array}{cc}
a_{\frac{N}{2}}\\
a_{\frac{N}{2}}
\end{array}
\right)
\bigg]\nn
\end{eqnarray}
\begin{eqnarray}
E_{-\frac{N}{2}}=\oplus_{n=1}^{N}\bigg[
(-1)^{n}\left(
\begin{array}{cc}
-b_{\frac{N}{2}}\\
b_{\frac{N}{2}}
\end{array}
\right)\bigg]\nn
\end{eqnarray}
are the eigensubspaces associated to the eigenvalues $\lambda_{\frac{N}{2}}= \frac{2}{h}$, $\lambda_{-\frac{N}{2}}=-\frac{2}{h}$ respectively, where $a_{\frac{N}{2}}, b_{\frac{N}{2}} \in {\bf R}^{m}, h=\frac{1}{N}$.

Since
\begin{eqnarray}
\mathcal{A}_{h,1}|_{E_{0}}=-{\hat B}_{j}=-
\left(
\begin{array}{cc}
Z_{j}&0\\
0&Z_{j}\\
\end{array}
\right),\nn
\end{eqnarray}
 whose eigenvalue $\Lambda_{1}=-\pi<0$ with multiplicity $2(m-1)$, and $\Lambda_{2}=-(j-m+1)\pi$ with multiplicity $2$. Note that from $j\ne m-1$, there hold $\Lambda_{2}>0$ when $j<m-1$, $\Lambda_{2}<0$ when $j>m-1$. We discuss the eigenvalues of $\mathcal{A}_{h}$ in two cases.

{\bf i).} $N$ is odd integer.

{\bf Subcase 1.1.1.} $j>m-1$.

For any $k\in \{1,2,\dots,[\frac{N}{2}]\}$, $E_{k}\oplus E_{-k}$ is $4m$-dimension invariant subspace of $\mathcal{A}_{h,1}$, and we have
\begin{eqnarray}
\mathcal{A}_{h,1}|_{E_{k}\oplus E_{-k}}=
\left(
\begin{array}{cc}
-\lambda _{k}I_{2m}-{\hat B}_{j}&0\\
0&-\lambda _{-k}I_{2m}-{\hat B}_{j}\\
\end{array}
\right).\nn
\end{eqnarray}
According to ${\hat B}_{j}$, we obtain that $\Lambda_{3}=2\frac{sin\theta_{k}}{h}-\pi>0$ and $\Lambda_{4}=-2\frac{sin\theta_{k}}{h}-\pi<0$ with multiplicity $2(m-1)$, $\Lambda_{5}=-2\frac{sin\theta_{k}}{h}-(j-m+1)\pi<0$ and $\Lambda_{6}=2\frac{sin\theta_{k}}{h}-(j-m+1)\pi$ with multiplicity $2$. Moreover since $N$ is sufficiently large, when $k<\frac{j-m+1}{2}$, $\Lambda_{6}<0$. When $k>\frac{j-m+1}{2}$, $\Lambda_{6}>0$. Summarize the above discussions, it follows that
\begin{eqnarray}
m_{j}^{+}&=&2(m-1)[\frac{N}{2}]+2([\frac{N}{2}]-\frac{j-m}{2})\nn\\
         &=&2(m-1)\frac{N-1}{2}+2(\frac{N-1}{2}-\frac{j-m}{2})\nn\\
         &=&(m-1)(N-1)+N-1-(j-m)\nn\\
         &=&mN-j.\nn
\end{eqnarray}

{\bf Subcase 1.1.2.} $j<m-1$.

Similar to subcase 1.1.1, we also have
$\Lambda_{3}=2\frac{sin\theta_{k}}{h}-\pi>0$ and $\Lambda_{4}=-2\frac{sin\theta_{k}}{h}-\pi<0$ with multiplicity $2(m-1)$, $\Lambda_{5}=-2\frac{sin\theta_{k}}{h}-(j-m+1)\pi$ and $\Lambda_{6}=2\frac{sin\theta_{k}}{h}-(j-m+1)\pi>0$ with multiplicity $2$. Moreover since $N$ is sufficiently large, when $k<\frac{m-j-1}{2}$, $\Lambda_{5}>0$. When $k>\frac{m-j-1}{2}$, $\Lambda_{5}<0$. Thus there hold
\begin{eqnarray}
m_{j}^{+}&=&2(m-1)[\frac{N}{2}]+2(\frac{m-j}{2}-1)+2+2[\frac{N}{2}]\nn\\
         &=&2(m-1)\frac{N-1}{2}+2(\frac{m-j}{2}-1)+2+2\frac{N-1}{2}\nn\\
         &=&(m-1)(N-1)+(m-j)-2+2+N-1\nn\\
         &=&mN-j.\nn
\end{eqnarray}

{\bf ii).} $N$ is even integer.

{\bf Subcase 1.1.3.} $j>m-1$.

For any $k\in \{1,2,\dots,\frac{N}{2}-1\}$, $E_{k}\oplus E_{-k}$ is $4m$-dimension invariant subspace of $\mathcal{A}_{h,1}$, and we have
\begin{eqnarray}
\mathcal{A}_{h,1}|_{E_{k}\oplus E_{-k}}=
\left(
\begin{array}{cc}
-\lambda _{k}I_{2m}-{\hat B}_{j}&0\\
0&-\lambda _{-k}I_{2m}-{\hat B}_{j}\\
\end{array}
\right).\nn
\end{eqnarray}
According to ${\hat B}_{j}$, we obtain that $\Lambda_{3}=2\frac{sin\theta_{k}}{h}-\pi>0$ and $\Lambda_{4}=-2\frac{sin\theta_{k}}{h}-\pi<0$ with multiplicity $2(m-1)$, $\Lambda_{5}=-2\frac{sin\theta_{k}}{h}-(j-m+1)\pi<0$ and $\Lambda_{6}=2\frac{sin\theta_{k}}{h}-(j-m+1)\pi$ with multiplicity $2$. Moreover since $N$ is sufficiently large, when $k<\frac{j-m+1}{2}$, $\Lambda_{6}<0$. when $k>\frac{j-m+1}{2}$, $\Lambda_{6}>0$.

For $k=\frac{N}{2}$, $E_{k}\oplus E_{-k}$ is $2m$-dimension invariant subspace of $\mathcal{A}_{h,1}$, and we have
\begin{eqnarray}
\mathcal{A}_{h,1}|_{E_{\frac{N}{2}}\oplus E_{-\frac{N}{2}}}=
\left(
\begin{array}{cc}
-\lambda _{k}I_{m}-Z_{j}&0\\
0&-\lambda _{-k}I_{m}-Z_{j}\\
\end{array}
\right).\nn
\end{eqnarray}
According to $Z_{j}$, we obtain that $\Lambda_{7}=\frac{2}{h}-\pi>0$ and $\Lambda_{8}=-\frac{2}{h}-\pi<0$ with multiplicity $m$. Summarize the above discussions, it follows that
\begin{eqnarray}
m_{j}^{+}&=&2(m-1)([\frac{N}{2}]-1)+m+2([\frac{N}{2}]-1-\frac{j-m}{2})\nn\\
         &=&2(m-1)(\frac{N}{2}-1)+m+2(\frac{N}{2}-1-\frac{j-m}{2})\nn\\
         &=&(m-1)(N-2)+m+(N-2-j+m)\nn\\
         &=&mN-j.\nn
\end{eqnarray}

{\bf Subcase 1.1.4.} $j<m-1$.

Similar to subcase 1.1.3, for any $k\in \{1,2,\dots,\frac{N}{2}-1\}$, we also have
$\Lambda_{3}=2\frac{sin\theta_{k}}{h}-\pi>0$ and $\Lambda_{4}=-2\frac{sin\theta_{k}}{h}-\pi<0$ with multiplicity $2(m-1)$, $\Lambda_{5}=-2\frac{sin\theta_{k}}{h}-(j-m+1)\pi$ and $\Lambda_{6}=2\frac{sin\theta_{k}}{h}-(j-m+1)\pi>0$ with multiplicity $2$. Moreover since $N$ is sufficiently large, when $k<\frac{m-j-1}{2}$, $\Lambda_{5}>0$. When $k>\frac{m-j-1}{2}$, $\Lambda_{5}<0$. In addition to $\Lambda_{7}=\frac{2}{h}-\pi>0$ and $\Lambda_{8}=-\frac{2}{h}-\pi<0$ with multiplicity $m$, thus there hold
\begin{eqnarray}
m_{j}^{+}&=&2+2(m-1)([\frac{N}{2}]-1)+2([\frac{N}{2}]-1)+2(\frac{m-j}{2}-1)+m\nn\\
         &=&2+2(m-1)(\frac{N}{2}-1)+2(\frac{N}{2}-1)+2(\frac{m-j}{2}-1)+m\nn\\
         &=&2+(m-1)(N-2)+(N-2)+m-j-2+m\nn\\
         &=&mN-j.\nn
\end{eqnarray}

{\bf Case 1.2.} $m=1$ and $j\in 2{\bf Z}\setminus \{0\}$.\\
In this case, $\hat {B}_{j}(t)$ is not a constant matrix. Since the proof is similar to Lemma 6.1.7 in \cite{Lon4}, we sketch it here. We can define $\gamma\in {\mathcal {P}^{*}_{1}}(4)$ which corresponding to the fundamental solution to the continuous Hamiltonian system (\ref{1.1}) with coefficient matrix
\bea
B(t)=\hat {B}_{j}(t)\diamond (\pi I_{2})=[-\dot {\omega}_{1}(t)j\pi I_{2}+(\dot {\omega}_{2}(t)ln2) K]\diamond (\pi I_{2}) \lb{1.55}
\eea
for $t\in [0,1]$ and satisfying that $i(\gamma)=j+1$, where $K$ is defined by (\ref{1.109}). According to (\ref{1.55}), Choose
\bea
B_{n}=\hat {B}_{j}(\frac{n}{N})\diamond (\pi I_{2})=[-\dot {\omega}_{1}(\frac{n}{N})j\pi I_{2}+(\dot {\omega}_{2}(\frac{n}{N})ln2) K]\diamond (\pi I_{2}) \lb{1.56}
\eea
as the coefficient matrix of the discrete Hamiltonian system (\ref{1.4}), denote by $\gamma_{d,h,N}$ the corresponding fundamental solution. Since $\nu(\gamma)=0$, then there exists a neighborhood of $\gamma(1)$ in $Sp^{*}_{1}(2m)$, hence there is a positive distance between $\gamma(1)$ and $Sp^{0}(2m)$. On the other hand, if  $N$ is large enough, then $\gamma_{d,h,N}(1)$ is sufficiently close to $\gamma(1)$, together with $\nu(\gamma_{d,h,N}(1))=0$, it is followed that $i(\gamma_{d,h,N})=i(\gamma)=j+1$. Denote by $m^{+}$ positive Morse index of the discrete Hamiltonian system whose coefficient matrix is the left hand side of (\ref{1.56}), by $m_{j}^{+}, m_{1}^{+}$ positive Morse index of the discrete Hamiltonian system whose coefficient matrix is the first or second term  right hand side of (\ref{1.56}) respectively. Thus we have
\bea
m^{+}=m_{j}^{+}+m_{1}^{+},\nn
\eea
then apply Case 1.1 of this theorem to yield that
\bea
2N-(j+1)=m_{j}^{+}+N-1,\nn
\eea
hence $m_{j}^{+}=N-j$, the result holds.

{\bf Case 2.} $\omega \in {\bf U}\setminus{\bf R}$.

 Since $\omega \in {\bf U}\setminus{\bf R}$, we have $\alpha \in (0,\pi)\cup (\pi,2\pi)$, $\hat{B}_{j}$ can be chosen as $(j\pi) I_{2}\diamond O_{2m-2}$, where $O_{2m-2}$ is $(2m-2)\times (2m-2)$ zero matrix. Then $\Lambda_{1}=2\frac{sin\alpha_{k}}{h}>0$ and $\Lambda_{2}=-2\frac{sin\alpha_{k}}{h}<0$ with multiplicity $(m-1)$, $\Lambda_{3}=-2\frac{sin\alpha_{k}}{h}-j\pi<0$ and $\Lambda_{4}=2\frac{sin\alpha_{k}}{h}-j\pi$ with multiplicity $1$, where $0 \le k \le N-1$.
We only discuss the case $j>0$, the case $j\le 0$  is similarly dealt with. The proof is divided into two subcases.

{\bf Subcase 2.1.} $\alpha \in (0,\pi)$.

i). $j$ is even. When $\frac{j}{2} \le k \le N-1-\frac{j}{2}$, we have $\Lambda_{4}>0$. hence $m_{j}^{+}=(m-1)N+N-j=mN-j$.

ii). $j$ is odd. When $\frac{j-1}{2}+1 \le k \le N-1-\frac{j-1}{2}$, we have $\Lambda_{4}>0$. hence $m_{j}^{+}=(m-1)N+N-(\frac{j-1}{2}+1+\frac{j-1}{2})=mN-j$.

{\bf Subcase 2.2}. $\alpha \in (\pi,2\pi)$.

i). $j$ is even. When $\frac{j}{2}-1 \le k \le N-2-\frac{j}{2}$, we have $\Lambda_{4}>0$. hence $m_{j}^{+}=(m-1)N+N-(\frac{j}{2}-1+\frac{j}{2}+1)=mN-j$.

ii). $j$ is odd. When $\frac{j-1}{2} \le k \le N-2-\frac{j-1}{2}$, we have $\Lambda_{4}>0$. hence $m_{j}^{+}=(m-1)N+N-(\frac{j-1}{2}+\frac{j-1}{2}+1)=mN-j$.

{\bf Case 3.} $\omega=-1$.

Since $\omega=-1$, we have $\alpha=\pi$. We only discuss the case $j>0$ , the case $j\le0$ is similarly dealt with.

{\bf Subcase 3.1.} $j$ is even. \\
$\hat{B}_{j}$ can be $(j\pi) I_{2}\diamond O_{2m-2}$, then $\Lambda_{1}=2\frac{sin\alpha_{k}}{h}>0$ and $\Lambda_{2}=-2\frac{sin\alpha_{k}}{h}<0$ with multiplicity $(m-1)$, $\Lambda_{3}=-2\frac{sin\alpha_{k}}{h}-j\pi<0$ and $\Lambda_{4}=2\frac{sin\alpha_{k}}{h}-j\pi$ with multiplicity $1$, where $0 \le k \le N-1$. When $\frac{j}{2} \le k \le N-1-\frac{j}{2}$, we have $\Lambda_{4}>0$. hence $m_{j}^{+}=(m-1)N+N-j=mN-j$.

{\bf Subcase 3.2.} $j$ is odd. \\
i). If $m\geq 2$, then $\hat{B}_{j}$ can be $(j\pi-\pi) I_{2}\diamond ({\pi} I_{2})\diamond O_{2m-4}$, then $\Lambda_{1}=2\frac{sin\alpha_{k}}{h}>0$ and $\Lambda_{2}=-2\frac{sin\alpha_{k}}{h}<0$ with multiplicity $(m-2)$, $\Lambda_{3}=-2\frac{sin\alpha_{k}}{h}-(j-1)\pi<0$ and $\Lambda_{4}=2\frac{sin\alpha_{k}}{h}-(j-1)\pi$ with multiplicity $1$, $\Lambda_{5}=-2\frac{sin\alpha_{k}}{h}-\pi<0$ and $\Lambda_{6}=2\frac{sin\alpha_{k}}{h}-\pi$ with multiplicity $1$, where $0 \le k \le N-1$. When $\frac{j-1}{2} \le k \le N-1-\frac{j-1}{2}$, we have $\Lambda_{4}>0$. Moreover $\Lambda_{6}>0$ for $k \in \{1, \dots, N-1\}$, hence $m_{j}^{+}=(m-2)N +(N-j+1)+N-1=mN-j$.\\
ii). If $m=1$, then $B(t)$ can be chosen as (\ref{1.55}), by a similar argument as Subcase 1.2, together with case 1 and i) of subcase 3.2, we obtain the result.

{\bf Step 2.} $\nu_{\omega}(\gamma_{d,h,N}(1))\ne 0.$\\
The problem is reduced to prove for the path $\gamma_{d,s}$, where $\gamma_{d,s}$ is given by (\ref{1.45}), note that $\gamma_{d,s}(1) \in Sp(2n)_{\omega}^*$ for $s\in [-1,1]\setminus\{0\}$. As case $\xi(r,\cdot)$ with $0 \le r \le r_{1}$ of Step 1, we can obtain the corresponding $B_{s,n}$ and discrete Hamiltonian system
\be
\frac{z_{n+1}-z_{n}}{h}=JB_{s,n}\tilde {z}_{n},\lb{1.147}
\ee
where $s \in [-1,1], 0 \le n \le N-1$.

By Lemma 2.9, we have
\bea
m^{0}_{h,\omega}=\nu_{\omega}(\gamma_{d,h,N}(1)).\lb{1.57}
\eea
Note that $\mathcal {A}_{s,h,\omega}$ is a small perturbation of $\mathcal {A}_{h,\omega}$, which are defined by (\ref{1.47}) and  (\ref{1.48}) respectively, thus we have
\bea
m^{+}_{h,\omega}\le m_{s,h,\omega}^{+},&& m^{-}_{h,\omega}\le m_{-s,h,\omega}^{-}.\lb{1.59}
\eea
Using the above theorem, for $s\in [-1,1]\setminus \{0\}$, we have
\bea
m_{s,h,\omega}^{-}=mN+i_{\omega}(\gamma_{d,s}),\hspace*{2mm} m_{s,h,\omega}^{0}=0, \hspace*{2mm}m_{s,h,\omega}^{+}=mN-i_{\omega}(\gamma_{d,s}).\lb{1.58}
\eea
For sufficiently small $s \in (0,1]$, combining (\ref{1.57}),(\ref{1.59}), (\ref{1.58}) with Lemma 4.2, we obtain
\begin{eqnarray}
mN-i_{\omega}(\gamma_{d,-s})&=& 2mN-m^{-}_{s,h,\omega} \le 2mN-m^{-}_{h,\omega}\nn\\
                            &=& m^{+}_{h,\omega}+m^{0}_{h,\omega}=m^{+}_{h,\omega}+\nu_{\omega}(\gamma_{d,h,N})\nn\\
                       &\le& mN-i_{\omega}(\gamma_{d,s})+i_{\omega}(\gamma_{d,s})-i_{\omega}(\gamma_{d,-s})\nn\\
                       &=& mN-i_{\omega}(\gamma_{d,-s})\lb{1.60}
\end{eqnarray}
therefore, it yields $m^{+}_{h,\omega}=mN-i_{\omega}(\gamma_{d,-s})-\nu_{\omega}(\gamma_{d,h,N})$, the proof is complete.\hfill\hb\\

{\bf Proof of Corollary 1.2} i). The result is followed from Theorem 1.1.\\
ii). Since $B_{n}=B(\frac{n}{N})$, the discrete system (\ref{1.66}) is the discretization of the continuous system (\ref{1.1}).  We conclude that $\gamma_{d,h,N}(1)$ and $\gamma(1)$ are in the same connected component of $Sp(2n)_{\omega}^*$ from $\nu_{\omega}(\gamma)=0$ and the appendix, thus $\nu_{\omega}(\gamma_{d,h,N})=\nu_{\omega}(\gamma)=0$ and $i_{\omega}(\gamma_{d,h,N})=i_{\omega}(\gamma)$, finally $Sign \mathcal{A}_{h,\omega}=2i_{\omega}(\gamma)$ holds from i).\hfill\hb

{\bf Proof of Corollary 1.3} i). By Theorem 1.1 and Definition 2.6, there exists
\begin{eqnarray}
\mathcal{S}_{h,\omega}^{\pm}&=&\lim _{\theta \rightarrow 0^{\pm}}m_{h,{\pm\theta}\omega e^{\sqrt -1}}^{-}-m_{h,\omega}^{-}\nn\\
          &=&\lim _{\theta \rightarrow 0^{\pm}}(mN+i_{{\pm\theta}\omega e^{\sqrt -1}}(\gamma_{d,h,N}))- (mN+i_{\omega}(\gamma_{d,h,N}))\nn\\
          &=&\lim _{\theta \rightarrow 0^{\pm}}(mN+i_{{\pm\theta}\omega e^{\sqrt -1}}(\gamma_{c,h,N}))- (mN+i_{\omega}(\gamma_{c,h,N}))\nn\\
          &=&\lim _{\theta \rightarrow 0^{\pm}}i_{{\pm\theta}\omega e^{\sqrt -1}}(\gamma_{c,h,N})- i_{\omega}(\gamma_{c,h,N})\nn\\
          &=&S_{\gamma_{c,h,N}(1)}^{\pm}(\omega),\nn
\end{eqnarray}
the second equality is followed from Theorem 1.1, the third equality uses the Definition 2.6, the last line holds from Definition 4.5. \\
ii). It follows that $i_{\omega}(\gamma_{d,h,N})=i_{\omega}(\gamma)$ from ii) of Corollary 1.2, together with i), the proof is complete.\hfill\hb

\setcounter{equation}{0}
\section{Appendix. The estimates of the errors}
For continuous system
\bea
{\dot z}(t)=f(t,z),\lb{1.111}
\eea
where $z(t)=(x^{T}(t),y^{T}(t))^{T}$ and $x(t), y(t) \in {\bf C}^{m},$ $t_{0} \le t \le T_{0}. $ The corresponding discrete system
\bea
\frac{z_{n+1}-z_{n}}{h}=f(t_{n},{\tilde z}_{n}),\lb{1.112}
\eea
where $0 < h=\frac{T_{0}-t_{0}}{N}$ with $N \in {\bf N}$ large enough, $z_{n}=(x^{T}_{n},y^{T}_{n})^{T}$, $\tilde {z}_{n}=(x^{T}_{n+1},y^{T}_{n})^{T}$ with $x_{n}, y_{n} \in {\bf C}^{m}$. Write $t_{n}=t_{0}+nh,$ for $0 \le n \le N$, moreover suppose that the system (\ref{1.111}) and (\ref{1.112}) possess the same initial value $ z(t_{0})=z_{0}$ and $f(\cdot, z)$ is continuous in $t$ and Lipschitz continuous in $z$ with Lipschitz constant $L$, i.e.
\bea
|f(t,z)-f(t,z')|\le L|z-z'|,~~~ \forall t \in [t_{0},T_{0}], \forall z, z' \in {\bf C}^{2m}.\lb{1.118}\eea
Set
\bea
\epsilon_{n}=z(t_{n})-z_{n}.\lb{1.116}
\eea

Moreover, we consider the continuous system
\bea
{\dot z}_{s}(t)=JB_{s}(t)z_{s}(t),\lb{1.170}
\eea
where $\dot z_{s}$ denote the derivative with respect to $t$, $z_{s}(t)=(x_{s}^{T}(t),y_{s}^{T}(t))^{T}$ is continuous in $s$, $x_{s}(t), y_{s}(t) \in {\bf C}^{m},$  $B \in C([0,1]\times[t_{0},T_{0}], \mathcal {L}_{s}({\bf R})^{2m}))$ with $\mathcal {L}_{s}({\bf R})^{2m}$ being the set of $2m \times 2m $ symmetric matrices, $t_{0} \le t \le T_{0}.$ The corresponding discrete system
\bea
\frac{z_{s,n+1}-z_{s,n}}{h}=JB_{s}(t_{n}){\tilde z}_{s,n},\lb{1.171}
\eea
where $0 < h=\frac{T_{0}-t_{0}}{N}$ with $N \in {\bf N}$ large enough, $z_{s,n}=(x^{T}_{s,n},y^{T}_{s,n})^{T}$, $\tilde {z}_{s,n}=(x^{T}_{s,n+1},y^{T}_{s,n})^{T}$ with $x_{s,n}, y_{s,n} \in {\bf C}^{m}$ and $t_{n}=t_{0}+nh$ for $0 \le n \le N$. Also suppose that the system (\ref{1.170}) and (\ref{1.171}) possess the same initial value $ z_{s}(t_{0})=z_{s,0}$.

Set
\bea
\epsilon_{s,n}=z_{s}(t_{n})-z_{s,n},\lb{1.172}
\eea
then we have

{\it{\bf A.1.}{~~\bf i).} $\epsilon_{n}=o(1)$ for $0 \le n \le N$ as $h\rightarrow 0$.

              {~~~~~~~\bf ii).} $\epsilon_{s,n}=\epsilon_{n}=o(1)$ for $0 \le n \le N$ independent on $s$ as $h\rightarrow 0$.

              {~~~~~~~\bf iii).} For $0 \le n \le N$, $z_{s,n}$ is continuous in $s$.}\\
{\bf Proof.} The proof of i) is divided into three steps.

{\bf Step 1.} From (\ref{1.111}), there holds
\begin{eqnarray}
z(t_{n+1})&=&z(t_{n})+\int_{t_{n}}^{t_{n}+h}f(\tau,z(\tau))d\tau \nn\\
          &=&z(t_{n})+hf(t_{n},z(t_{n}))+R_{n}(h),\lb{1.113}
\end{eqnarray}
where
\bea
R_{n}(h)=\int_{t_{n}}^{t_{n}+h}f(\tau,z(\tau))d\tau-hf(t_{n},z(t_{n})),\lb{1.114}
\eea
Note that for $0 \le n \le N$ there holds $R_{n}(h)=o(h)$.

Meanwhile, by (\ref{1.112}), we have
\bea
z_{n+1}=z_{n}+hf(t_{n},{\tilde z}_{n}).\lb{1.115}
\eea

{\bf Step 2.} (\ref{1.113}) is subtracted by (\ref{1.115}), together with the definition (\ref{1.116}) of $\epsilon_{n}$, we get that
\bea
\epsilon_{n+1}=\epsilon_{n}+hf(t_{n},z(t_{n}))- hf(t_{n},{\tilde z}_{n})+R_{n}(h).\lb{1.117}
\eea
Firstly, since
\begin{eqnarray}
|x(t_{n})-x_{n+1}|&=&|x(t_{n})-x(t_{n+1})+x(t_{n+1})-x_{n+1}|\nn\\
                  &\le&|x(t_{n})-x(t_{n+1})|+|x(t_{n+1})-x_{n+1}|\nn\\
                  &=&|x(t_{n})-x(t_{n+1})|+|\epsilon_{n+1}|\nn\\
                  &\le&Mh+|\epsilon_{n+1}|.\lb{1.120}
\end{eqnarray}
where the last line holds from the fact ${\dot z}(t)$ is continuous in $t$, hence we can set $M=max_{t_{0}\le t \le T_{0}}|{\dot z}(t)|$. Therefore we can estimate the second term and third term of the right hand side of (\ref{1.117})
\begin{eqnarray}
&&|f(t_{n},z(t_{n}))-f(t_{n},{\tilde z}_{n})|\nn\\
&&\le L|z(t_{n})-{\tilde z}_{n}|\nn\\
&&=L\bigg|\left(
\begin{array}{cc}
x(t_{n})-x_{n+1}\\
y(t_{n})-y_{n}\\
\end{array}
\right)\bigg|\nn\\
&&\le L(|\epsilon_{n}|+|\epsilon_{n+1}|+hM),\lb{1.119}
\end{eqnarray}
where we use (\ref{1.118}) in the first inequality, the third inequality is followed from (\ref{1.120}).

{\bf Step 3.} Since the finiteness of $N$, there exists certain function $R(h)=o(h)$ satisfying $|R_{n}(h)| \le R(h)$ for $0 \le n \le N$. For example, we can choose $R(h)=max_{0\le n \le N}R_{n}(h)$. Hence from (\ref{1.117}) and (\ref{1.119}), we have
\bea
|\epsilon_{n+1}|&\le&|\epsilon_{n}|+hL(|\epsilon_{n}|+|\epsilon_{n+1}|+hM)+|R_{n}(h)|\nn\\
                &\le & \epsilon_{n}|+hL|\epsilon_{n}|+hL|\epsilon_{n+1}|+h^{2}LM)+|R(h)|.\lb{1.121}
\eea
This gives
\bea
(1-hL)|\epsilon_{n+1}|\le (1+hL)|\epsilon_{n}|+h^{2}LM+R(h)\nn
\eea
 Since $h=\frac{1}{N}$ is small enough, then $1-hL>0$, thus we see that
\bea
|\epsilon_{n+1}|&\le&\frac{1+hL}{1-hL}|\epsilon_{n}|+\frac{h^{2}LM}{1-hL}+\frac{R(h)}{1-hL}\nn\\
                &\le&\frac{1+hL}{1-hL}\bigg[\frac{1+hL}{1-hL}|\epsilon_{n}|+\frac{h^{2}LM}{1-hL}+\frac{R(h)}{1-hL}\bigg]+\frac{h^{2}LM}{1-hL}+\frac{R(h)}{1-hL}\nn\\
                &=&(\frac{1+hL}{1-hL})^2|\epsilon_{n-1}|+\frac{h^{2}LM}{1-hL}\bigg[\frac{1+hL}{1-hL}+1\bigg]+\frac{R(h)}{1-hL}\bigg[\frac{1+hL}{1-hL}+1\bigg].\lb{1.122}
\eea
By an induction argument, it forces from (\ref{1.122}) that
\bea
|\epsilon_{n+1}|\le(\frac{1+hL}{1-hL})^{n+1}|\epsilon_{0}|+\frac{h^{2}LM}{1-hL}\Sigma_{j=0}^{n}(\frac{1+hL}{1-hL})^{j}+\frac{R(h)}{1-hL}\Sigma_{j=0}^{n}(\frac{1+hL}{1-hL})^{j}.\lb{1.123}
\eea
Replace $n+1$ by $n$ in (\ref{1.123}), and note that the fact that $t_{n}=t_{0}+nh \le T_{0}$ yields $n\le \frac{T_{0}-t_{0}}{h}$, there holds
\bea
|\epsilon_{n}|&\le&(\frac{1+hL}{1-hL})^{n}|\epsilon_{0}|+\frac{h^{2}LM}{1-hL}\Sigma_{j=0}^{n-1}(\frac{1+hL}{1-hL})^{j}+\frac{R(h)}{1-hL}\Sigma_{j=0}^{n-1}(\frac{1+hL}{1-hL})^{j}\nn\\
              &\le&(\frac{1+hL}{1-hL})^{\frac{T_{0}-t_{0}}{h}}|\epsilon_{0}|+\frac{h^{2}LM}{1-hL}\frac{(\frac{1+hL}{1-hL})^{n}-1}{\frac{1+hL}{1-hL}-1}+\frac{R(h)}{1-hL}\frac{(\frac{1+hL}{1-hL})^{n}-1}{\frac{1+hL}{1-hL}-1}\nn\\
              &\le& 2e^{2\frac{T_{0}-t_{0}}{1-hL}}|\epsilon_{0}|+\frac{h^{2}LM}{1-hL}\frac{1-hL}{hL}\bigg[e^{2\frac{T_{0}-t_{0}}{1-hL}}-1\bigg] +\frac{R(h)}{1-hL}\frac{1-hL}{hL}\bigg[e^{2\frac{T_{0}-t_{0}}{1-hL}}-1\bigg]\nn\\
              &\le& 2e^{2\frac{T_{0}-t_{0}}{1-hL}}|\epsilon_{0}|+hM\bigg[e^{2\frac{T_{0}-t_{0}}{1-hL}}-1\bigg] +\frac{R(h)}{hL}\bigg[e^{2\frac{T_{0}-t_{0}}{1-hL}}-1\bigg].\lb{1.124}
\eea
In addition to  $R(h)=o(h)$ and $\epsilon_{0}=0$, the result is followed from (\ref{1.124}).\\
{\bf ii).} The proof is also divided into three steps.

{\bf Step 1.}
Note that since $B \in C([0,1]\times[t_{0},T_{0}], \mathcal {L}_{s}({\bf R})^{2m}))$, there exists $L_{0}$ large enough such that
\bea
|JB_{s}(t)z-JB_{s}(t)z'|\le L_{0}|z-z'|,~~~ \forall s \in [0,1],~~\forall t \in [t_{0},T_{0}], ~~\forall z, z' \in {\bf C}^{2m}.\lb{1.148}
\eea
From (\ref{1.170}), there holds
\begin{eqnarray}
z_{s}(t_{n+1})&=&z_{s}(t_{n})+\int_{t_{n}}^{t_{n}+h}f(s,\tau,z(\tau))d\tau \nn\\
          &=&z_{s}(t_{n})+hf(s,t_{n},z(t_{n}))+R_{s,n}(h),\lb{1.173}
\end{eqnarray}
where
\bea
R_{s,n}(h)=\int_{t_{n}}^{t_{n}+h}f(s,\tau,z(\tau))d\tau-hf(s,t_{n},z(t_{n})),\lb{1.174}
\eea
and $f(s,t,z)=JB_{s}(t)z_{s}$. Note that for any $s \in [0,1]$ and $0 \le n \le N$, we have $ R_{s,n}(h)=o(h)$ and $R_{s,n}(h)$ is continuous in $s$, by the the compactness of $[0,1]$, we conclude that there exists some $R_{n}(h)$ such that for $0 \le n \le N$ there uniformly holds
\bea
R_{s,n}(h)=R_{n}(h)=o(h)\lb{1.176}
\eea
for $0\le s \le 1$.

Meanwhile, by (\ref{1.171}), we have
\bea
z_{s,n+1}=z_{s,n}+hf(s,t_{n},{\tilde z}_{n}).\lb{1.175}
\eea

{\bf Step 2.} (\ref{1.173}) is subtracted by (\ref{1.175}), together with the definition (\ref{1.172}) of $\epsilon_{s,n}$, we get that
\bea
\epsilon_{s,n+1}=\epsilon_{s,n}+hf(s,t_{n},z(t_{n}))- hf(s,t_{n},{\tilde z}_{n})+R_{s,n}(h).\lb{1.176}
\eea
Firstly, since
\begin{eqnarray}
|x_{s}(t_{n})-x_{s,n+1}|&=&|x_{s}(t_{n})-x_{s}(t_{n+1})+x_{s}(t_{n+1})-x_{s,n+1}|\nn\\
                  &\le&|x_{s}(t_{n})-x_{s}(t_{n+1})|+|x_{s}(t_{n+1})-x_{s,n+1}|\nn\\
                  &=&|x_{s}(t_{n})-x_{s}(t_{n+1})|+|\epsilon_{s,n+1}|\nn\\
                  &\le&Mh+|\epsilon_{s,n+1}|.\lb{1.177}
\end{eqnarray}
where the last line holds from the fact ${\dot z}_{s}(t)$ is continuous in $s,t$ and the compactness of $[0,1]$, hence we can set $M=max_{t_{0}\le t \le T_{0}\atop 0 \le s \le 1}|{\dot z}_{s}(t)|$. Therefore we can estimate the second term and third term of the right hand side of (\ref{1.176})
\begin{eqnarray}
&&|f(s,t_{n},z(t_{n}))-f(s,t_{n},{\tilde z}_{n})|\nn\\
&&\le L_{0}|z_{s}(t_{n})-{\tilde z}_{s,n}|\nn\\
&&=L_{0}\bigg|\left(
\begin{array}{cc}
x_{s}(t_{n})-x_{s,n+1}\\
y_{s}(t_{n})-y_{s,n}\\
\end{array}
\right)\bigg|\nn\\
&&\le L_{0}(|\epsilon_{s,n}|+|\epsilon_{s,n+1}|+hM),\lb{1.178}
\end{eqnarray}
where we use (\ref{1.148}) in the first inequality, the third inequality is followed from (\ref{1.177}).

{\bf Step 3.} Since the finiteness of $N$, there exists certain function $R(h)=o(h)$ satisfying $|R_{n}(h)| \le R(h)$ for $0 \le n \le N$. For example, we can choose $R(h)=max_{0\le n \le N}R_{n}(h)$. Hence from (\ref{1.176}) and (\ref{1.178}), we have
\bea
|\epsilon_{s,n+1}|&\le&|\epsilon_{s,n}|+hL_{0}(|\epsilon_{s,n}|+|\epsilon_{s,n+1}|+hM)+|R_{n}(h)|\nn\\
                &\le & \epsilon_{s,n}|+hL_{0}|\epsilon_{s,n}|+hL_{0}|\epsilon_{s,n+1}|+h^{2}L_{0}M)+|R(h)|.\lb{1.179}
\eea
This gives
\bea
(1-hL_0)|\epsilon_{s,n+1}|\le (1+hL_0)|\epsilon_{s,n}|+h^{2}L_0M+R(h)\nn
\eea
Since $h=\frac{1}{N}$ is small enough, then $1-hL_0>0$, thus we see that
\bea
|\epsilon_{s,n+1}|&\le&\frac{1+hL_0}{1-hL_0}|\epsilon_{s,n}|+\frac{h^{2}L_0M}{1-hL_0}+\frac{R(h)}{1-hL_0}\nn\\
                &\le&\frac{1+hL_0}{1-hL_0}\bigg[\frac{1+hL_0}{1-hL_0}|\epsilon_{s,n}|+\frac{h^{2}L_0M}{1-hL_0}+\frac{R(h)}{1-hL_0}\bigg]+\frac{h^{2}L_0M}{1-hL_0}+\frac{R(h)}{1-hL_0}\nn\\
                &=&(\frac{1+hL_0}{1-hL_0})^2|\epsilon_{s,n-1}|+\frac{h^{2}L_0M}{1-hL_0}\bigg[\frac{1+hL_0}{1-hL_0}+1\bigg]+\frac{R(h)}{1-hL_0}\bigg[\frac{1+hL_0}{1-hL_0}+1\bigg].\lb{1.180}
\eea
By an induction argument, it forces from (\ref{1.180}) that
\bea
|\epsilon_{s,n+1}|\le(\frac{1+hL_0}{1-hL_0})^{n+1}|\epsilon_{s,0}|+\frac{h^{2}L_0M}{1-hL_0}\Sigma_{j=0}^{n}(\frac{1+hL_0}{1-hL_0})^{j}+\frac{R(h)}{1-hL_0}\Sigma_{j=0}^{n}(\frac{1+hL_0}{1-hL_0})^{j}.\lb{1.181}
\eea
Replace $n+1$ by $n$ in (\ref{1.181}), and note that the fact that $t_{n}=t_{0}+nh \le T_{0}$ yields $n\le \frac{T_{0}-t_{0}}{h}$, there holds
\bea
|\epsilon_{s,n}|&\le&(\frac{1+hL_0}{1-hL_0})^{n}|\epsilon_{s,0}|+\frac{h^{2}L_0M}{1-hL_0}\Sigma_{j=0}^{n-1}(\frac{1+hL_0}{1-hL_0})^{j}+\frac{R(h)}{1-hL_0}\Sigma_{j=0}^{n-1}(\frac{1+hL_0}{1-hL_0})^{j}\nn\\
              &\le&(\frac{1+hL_0}{1-hL_0})^{\frac{T_{0}-t_{0}}{h}}|\epsilon_{0}|+\frac{h^{2}L_0M}{1-hL_0}\frac{(\frac{1+hL_0}{1-hL_0})^{n}-1}{\frac{1+hL_0}{1-hL_0}-1}+\frac{R(h)}{1-hL_0}\frac{(\frac{1+hL_0}{1-hL_0})^{n}-1}{\frac{1+hL_0}{1-hL_0}-1}\nn\\
              &\le& 2e^{2\frac{T_{0}-t_{0}}{1-hL_0}}|\epsilon_{0}|+\frac{h^{2}L_0M}{1-hL_0}\frac{1-hL_0}{hL_0}\bigg[e^{2\frac{T_{0}-t_{0}}{1-hL_0}}-1\bigg] +\frac{R(h)}{1-hL_0}\frac{1-hL_0}{hL_0}\bigg[e^{2\frac{T_{0}-t_{0}}{1-hL_0}}-1\bigg]\nn\\
              &\le& 2e^{2\frac{T_{0}-t_{0}}{1-hL_0}}|\epsilon_{0}|+hM\bigg[e^{2\frac{T_{0}-t_{0}}{1-hL_0}}-1\bigg] +\frac{R(h)}{hL_0}\bigg[e^{2\frac{T_{0}-t_{0}}{1-hL_0}}-1\bigg].\lb{1.182}
\eea
In addition to  $R(h)=o(h)$ and $\epsilon_{0}=0$, the result is followed from (\ref{1.182}).\\

{\bf iii).} Since $z_{s}(t)$ is continuous in $s$, specifically $z_{s}(0)$ is continuous in $s$. While $z_{s,0}=z_{s}(0)$, thus $z_{s,0}$ is continuous in $s$, according to (\ref{1.171}), together with the fact that $B_{s,n}$ is continuous in $s$,  it yields that $z_{s,1}$ is continuous in $s$. By the induction argument, the result holds.  \hfill\hb\\
{\bf Acknowledgements.} I am deeply indebted to my thesis advisor, Professor Yiming Long for suggesting this
problem and for sharing his time and great expertise in index theory and Hamiltonian dynamics during numerous enlightening and stimulating discussions.

\bibliographystyle{abbrv}

\end{document}